\documentclass[10pt]{article}
\usepackage{amsmath,amssymb}
\usepackage[margin=1in]{geometry}
\usepackage{color}
\usepackage{graphicx}
\usepackage{url}
\usepackage{afterpage}
\newtheorem{lemma}{Lemma}
\newtheorem{remark}{Remark}
\newtheorem{definition}{Definition}
\def\ubar{\bar u}
\def\vbar{\bar v}
\def\qed{\hfill$\Box$}

\newcommand{\RR}{\mathbb{R}}

\newcommand{\paren}[1]{\left( #1 \right)}

\newcommand{\wt}[1]{\widetilde{#1}}
\newcommand{\wb}[1]{\overline{#1}}
\newcommand{\norm}[1]{\left\lVert#1\right\rVert}
\newcommand{\OO}{\mathcal{O}}

\newcommand{\proj}{\mathcal{P}}
\newcommand{\projK}{\mathcal{P}^K}
\newcommand{\syn}{\mathcal{Q}}
\newcommand{\synK}{\mathcal{Q}^K}
\newcommand{\calS}{\mathcal{S}}
\newcommand{\calB}{\mathcal{B}}
\newcommand{\calBK}{\mathcal{B}^K}

\begin{document}

\title{A fast solver for the narrow capture and narrow escape problems in
the sphere}

\author{Jason Kaye
  \thanks{Courant Institute of Mathematical Sciences,
    New York University, New York, New York 10012.
    (Email: {jkaye@cims.nyu.edu}).}
  \and
  Leslie Greengard
  \thanks{Courant Institute of Mathematical Sciences,
    New York University, New York, New York 10012;
    Flatiron Institute, Simons Foundation, New York, New York 10010.
    (Email: {greengard@cims.nyu.edu}).}
}
\date{}

\maketitle

\begin{abstract}

We present an efficient method to solve the narrow 
capture and narrow escape problems for the sphere.
The {\em narrow capture} problem models the equilibrium behavior of a Brownian
particle in the exterior of a sphere whose surface is reflective,
except for a collection of small absorbing patches.
The {\em narrow escape} problem is the dual problem: it models the behavior of a Brownian
particle confined to the interior of a sphere whose surface is reflective,
except for a collection of small patches through which it can escape.

Mathematically, these give rise to mixed Dirichlet/Neumann boundary
value problems of the Poisson equation.  They are numerically
challenging for two main reasons: (1) the solutions are non-smooth at
Dirichlet-Neumann interfaces, and (2) they involve adaptive mesh
refinement and the solution of large, ill-conditioned linear systems when 
the number of small patches is large.

By using the Neumann Green's functions for the sphere, 
we recast each boundary value problem as a system of first-kind integral
equations on the collection of patches. A block-diagonal 
preconditioner together with a multiple scattering formalism leads
to a well-conditioned system of second-kind integral equations and
a very efficient approach to discretization.
This system is solved iteratively using GMRES. We develop a hierarchical, 
fast multipole method-like algorithm to accelerate each matrix-vector product.
Our method is insensitive to the patch size, and the total cost scales
with the number $N$ of patches as $\mathcal{O}(N \log N)$, after a
precomputation whose cost depends only on the patch size and not on the
number or arrangement of patches. We demonstrate the method with several
numerical examples, and are able to achieve highly accurate solutions with
100\, 000 patches in one hour on a 60-core workstation. For that case,
adaptive discretization of each patch would lead to a dense linear
system with about 360 million degrees of freedom. Our preconditioned
system uses only 13.6 million ``compressed" degrees of freedom and a few dozen 
GMRES iterations.
\end{abstract}

\section{Introduction}

\begin{figure}[ht]
  \centering
    \includegraphics[width=.4\linewidth]{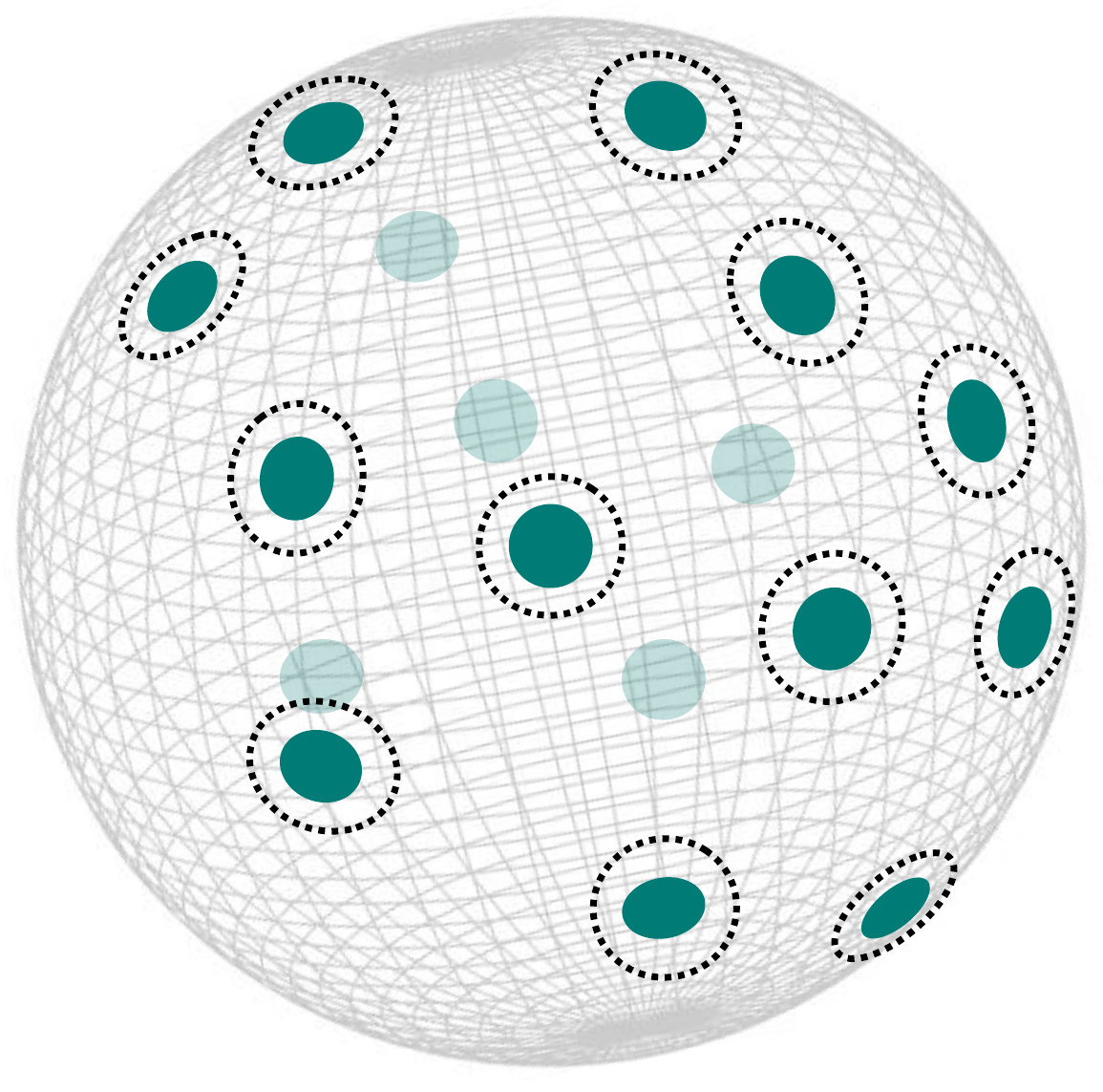}
    \caption{A sphere partially covered by disk-shaped patches. We
      assume each patch is of radius $\varepsilon$. We also assume that distinct
      patches are separated by a distance of at least $\varepsilon$.
      In the figure, this means that the regions bounded by the dashed
    lines do not overlap.}
    \label{fig:patches}
\end{figure}

We consider the numerical solution of two related
problems which arise in the study of Brownian diffusion by a particle in
the exterior or interior of a porous sphere. We denote the open unit
ball centered at the origin 
in $\mathbb{R}^3$ by $\Omega$, and assume that the sphere 
$\partial \Omega$ is
partially covered by $N$ small patches of radius $\varepsilon$, measured
in arclength (Fig. \ref{fig:patches}).
For the sake of simplicity, we assume that the patches are 
disk-shaped and comment briefly on more general shapes in the conclusion. 

The union of the patches is referred to as the {\em absorbing boundary}
and denoted by $\Gamma_A$.  The remainder of the boundary, 
$\Gamma_R = \partial \Omega \backslash \Gamma_A$, is referred to as the
{\em reflecting boundary}.
The first problem, called the narrow capture problem, is to calculate
the concentration $\ubar(x)$, at equilibrium, of Brownian particles at $x
\in \mathbb{R}^3 \backslash \wb{\Omega}$ with a given fixed concentration far
from the origin, assuming that particles are absorbed (removed) 
at $\Gamma_A$.  The second problem, 
called the narrow escape problem, is to
calculate the mean first passage time (MFPT) in $\Omega$, namely the
expected time $\vbar(x)$ for a Brownian particle released at $x \in \Omega$
to first reach $\Gamma_A$. In both settings, particles are 
reflected from $\Gamma_R$. In this paper, we sometimes refer to the
narrow capture problem as the exterior problem, and the narrow escape
problem as the interior problem.

These problems have received quite a lot of attention in the mathematics
and biophysics communities since the seminal work of Berg and Purcell
\cite{Berg1977}.  We do not seek to review the biophysical background
here, but note that the absorbing patches serve as a simplified model
for either surface receptors (the capture mechanism) or pores (the
escape mechanism) in an otherwise impermeable membrane. We refer the
reader to \cite{Berg1977,wiegel83,bressloff13,Holcmanbook2015,Holcman2014,kaizu14,isaacson11} for more detailed
discussions of applications and a selection of work on related biophysical
models.

Standard arguments from stochastic analysis show that both $\ubar$ and $\vbar$
satisfy a Poisson equation with mixed Dirichlet-Neumann boundary
conditions \cite{redner01,pavliotis14}. More precisely, for the capture
problem, if the far-field particle concentration is set to be $1$, 
then $\ubar$ satisfies the exterior Laplace equation:
\begin{align}\label{eq:extBVP}
  \begin{cases}
    \Delta \ubar = 0 & x \in \mathbb{R}^3 \backslash \wb{\Omega} \\
    \ubar = 0 & x \in \Gamma_A \\
    \frac{\partial \ubar}{\partial n} = 0 & x \in \Gamma_R \\
    \ubar(x) \to 1 & |x| \to \infty. \\
  \end{cases}
\end{align} 
A scalar quantity of interest is the
total equilibrium flux $J$ of particles through $\Gamma_A$:
\begin{equation}\label{eq:flux}
  J = \int_{\Gamma_A} \frac{\partial \ubar}{\partial n} \, dS.
\end{equation}
This is sometimes referred to as the capacitance of the system
(see Remark \ref{rem:capacitance}). 
For the escape problem, the MFPT $\vbar$ satisfies the interior Poisson equation:
\begin{align}\label{eq:intBVP}
  \begin{cases}
    \Delta \vbar = -1 & x \in \Omega \\
    \vbar = 0 & x \in \Gamma_A \\
    \frac{\partial \vbar}{\partial n} = 0 & x \in \Gamma_R.
  \end{cases}
\end{align}
Here, the quantity of interest is the average MFPT $\mu$ - that is the
average, over all possible initial particle positions, of the expected
time to escape from $\Omega$ through $\Gamma_A$:
\begin{equation}\label{eq:avgmfpt}
  \mu = \frac{1}{|\Omega|} \int_\Omega \vbar \, dV.
\end{equation}

Here, and in the remainder of the paper, 
$\frac{\partial}{\partial n}$ refers to the derivative in the outward
normal direction;
$n$ points towards the interior of $\Omega$ for the exterior problem, 
and towards the exterior of $\Omega$ for the interior problem.
In order to understand how the distribution of absorbing
patches on the surface affects $\ubar(x)$, $\vbar(x)$ and the associated
quantities $J$ and $\mu$,
a variety of asymptotic and numerical methods have been developed (see
\cite{Berg1977,Cheviakov2010,Lindsay2017,Bernoff2018,Singer2006,Holcmanbook2015,Holcman2014}
and the references therein). 

\begin{remark}\label{rem:capacitance}
  The total flux $J$ defined in \eqref{eq:flux} is sometimes referred
  to as the capacitance because of a connection to
  electrostatics. Imagine that
  the ball $\Omega$ is a dielectric with 
  low permittivity, and that $\Gamma_A$ is a collection of perfectly conducting
  patches on its surface, connected by infinitesimally thin wires so
  that they act as a single conductor. Suppose also that this object
  is surrounded by a dielectric with high permittivity and that the outer
  dielectric is enclosed by an infinitely large perfectly conducting sphere, 
  with a unit voltage drop from
  the outer conductor to the conducting patches. Then, letting the
  ratio of the permittivity of the outer dielectric to that of
  the inner dielectric approach $\infty$, the electrostatic potential
  outside $\wb{\Omega}$ satisfies \eqref{eq:extBVP}, and the
  electrostatic capacitance of the system is given by $J$.
\end{remark}

\begin{remark}\label{rem:mfpteqn}
  The total flux $J$ 
  is computed directly from the Neumann data on
  $\Gamma_A$, as seen from \eqref{eq:flux}. Likewise, the average MFPT $\mu$ can be computed
  directly from the Dirichlet data $\vbar$ on $\Gamma_R$. For this, we
  use Green's second identity,
  \[\int_\Omega \paren{\psi \Delta \varphi - \varphi \Delta
    \psi} \, dV = \int_{\partial \Omega} \paren{\psi
    \frac{\partial \varphi}{\partial n} - \varphi \frac{\partial
\psi}{\partial n}}\, dS\]
  with $\psi(x) \equiv \vbar(x)$ and $\varphi(x) \equiv \frac{|x|^2}{6}$.
  Using that $\Delta \frac{|x|^2}{6} = 1$, 
  $\int_\Omega \frac{|x|^2}{6} dV(x) =
  \frac{2 \pi}{15}$, and that for $|x| = 1$, $n
  \equiv x$ and $\frac{\partial}{\partial n} \frac{|x|^2}{6} = \frac13$, 
we obtain
  \[\int_\Omega \vbar \, dV = \frac13 \int_{\partial \Omega} \vbar \, dS - \frac16
    \int_{\partial \Omega} \frac{\partial \vbar}{\partial n} \, dS - \frac{2
    \pi}{15}.\]
  Applying the divergence theorem to the second term, dividing by
  $|\Omega|$, and using that $|\Omega| = \frac{4 \pi}{3}$,
  $|\partial \Omega| = 4 \pi$ gives an alternative expression for
  $\mu$:
  \begin{equation}\label{eq:avgmfpt2}
    \mu = \frac{1}{|\partial \Omega|}
    \int_{\partial \Omega} \vbar \, dS + \frac{1}{15} \equiv \frac{1}{|\partial \Omega|}
    \int_{\Gamma_R} \vbar \, dS + \frac{1}{15}.
  \end{equation}
  Thus the average MFPT over $\Omega$ may be obtained from the average
  MFPT on $\partial \Omega$.
\end{remark}

Given an arrangement of patches, we present here a
fast, high-order accurate numerical scheme for the evaluation
of $\ubar$, $J$, $\vbar$, and $\mu$, 
of particular use when $N$ is large and
$\varepsilon$ is small. 
Such computations are numerically challenging, partly because
solutions of elliptic boundary value problems of mixed type 
are singular near Dirichlet-Neumann interfaces
\cite{Sneddon1966,Fabrikant1989}.
Direct discretization, using either PDE-based methods or integral
equation methods, would require many degrees of freedom to resolve the 
singularities in $\ubar$ and $\vbar$. Further, the resulting linear systems would
be large and ill-conditioned, especially in cases involving large
numbers of small patches. 

The formulation presented here is
{\em well-conditioned}, is nearly identical for the capture and escape
problems, and suffers no loss in accuracy
or increase in computational cost as $\varepsilon$ is decreased. 
To make large-scale problems practical,
we have developed a fast algorithm, so that 
the cost per GMRES iteration \cite{GMRES} is of the order
$\OO(N \log N)$, rather than 
$\OO(N^2)$.
Our method involves the following ingredients:

\begin{itemize}
  \item We make use of the Neumann Green's functions for the 
    interior and exterior
    of the sphere to recast \eqref{eq:extBVP} and \eqref{eq:intBVP} as
    first-kind integral equations for a density $\sigma$ on
    $\Gamma_A$. 
  \item
    Given a patch radius $\varepsilon$, 
    we precompute the solution operator for the
    corresponding one-patch integral equation, assuming smooth Dirichlet
    data which is expanded in a rapidly converging series of
    Zernike polynomials. We analytically incorporate a square root singularity in the 
    induced density at the Dirichlet/Neumann interface. 
  \item To solve the many-patch integral equation, we use 
    the solution operator for the one-patch integral equation as 
    a block-diagonal ``right preconditioner''. This yields 
    a second-kind Fredholm system of equations which, upon discretization, is
    well-conditioned and has a small number of degrees of freedom per patch.
  \item We solve the resulting linear system by iteration, using GMRES, 
    and accelerate each matrix-vector product by means of a fast algorithm
    modeled after the fast multipole method (FMM). The fast algorithm uses the 
    {\em interpolative decomposition} \cite{liberty07} to derive a compressed
    representation of the outgoing field induced by the density on a patch, 
    a hierarchical organization of patches into groups at different length
    scales, and a spectral representation of the smooth incoming field due to
    densities on distant patches.
\end{itemize}

Though most of the past work on the narrow capture and narrow escape
problems is based on asymptotics, we wish to highlight the numerical work
of Bernoff and Lindsay, who also proposed an integral equation method
for the narrow capture problem for the sphere and the plane based on the
Neumann Green's function \cite{Bernoff2018}.  Our approach to discretization
shares several characteristics with theirs: both methods 
incorporate a square root singularity into the density on each patch 
analytically,
and both use a representation in terms of Zernike polynomials for
smooth Dirichlet data on each patch.

The paper is organized as follows. 
In Section \ref{sec:analyticalsetup}, we introduce the analytical 
framework for our method, reformulate the 
boundary value problems as first-kind integral
equations using single layer potentials, and explain 
how to calculate the scalar quantities $J$ and
$\mu$ directly as functionals of the layer potential densities. In
Section \ref{sec:manypatchsetuptools}, we show how to transform the
first-kind integral equations into Fredholm equations of the second-kind,
using the solution operator for the one-patch integral
equation as a preconditioner. 
In Sections \ref{sec:zernike}, \ref{sec:onepatchprelim}, and
\ref{sec:manypatchsetup} we describe our discretization approach for the
full system of equations, and in Section \ref{sec:outgoingrep} we
introduce the technical tools involved in our fast algorithm.
In Section \ref{sec:algorithm} we describe the full method, including
our fast algorithm to accelerate the application of the system matrix.
In Section \ref{sec:onepatch}, we provide a detailed description of the 
solver for the one-patch integral equation. We demonstrate the
performance of the method with numerical experiments
in Section \ref{sec:numresults}.

\begin{figure}[p!]
  \centering
    \includegraphics[width=.82\linewidth]{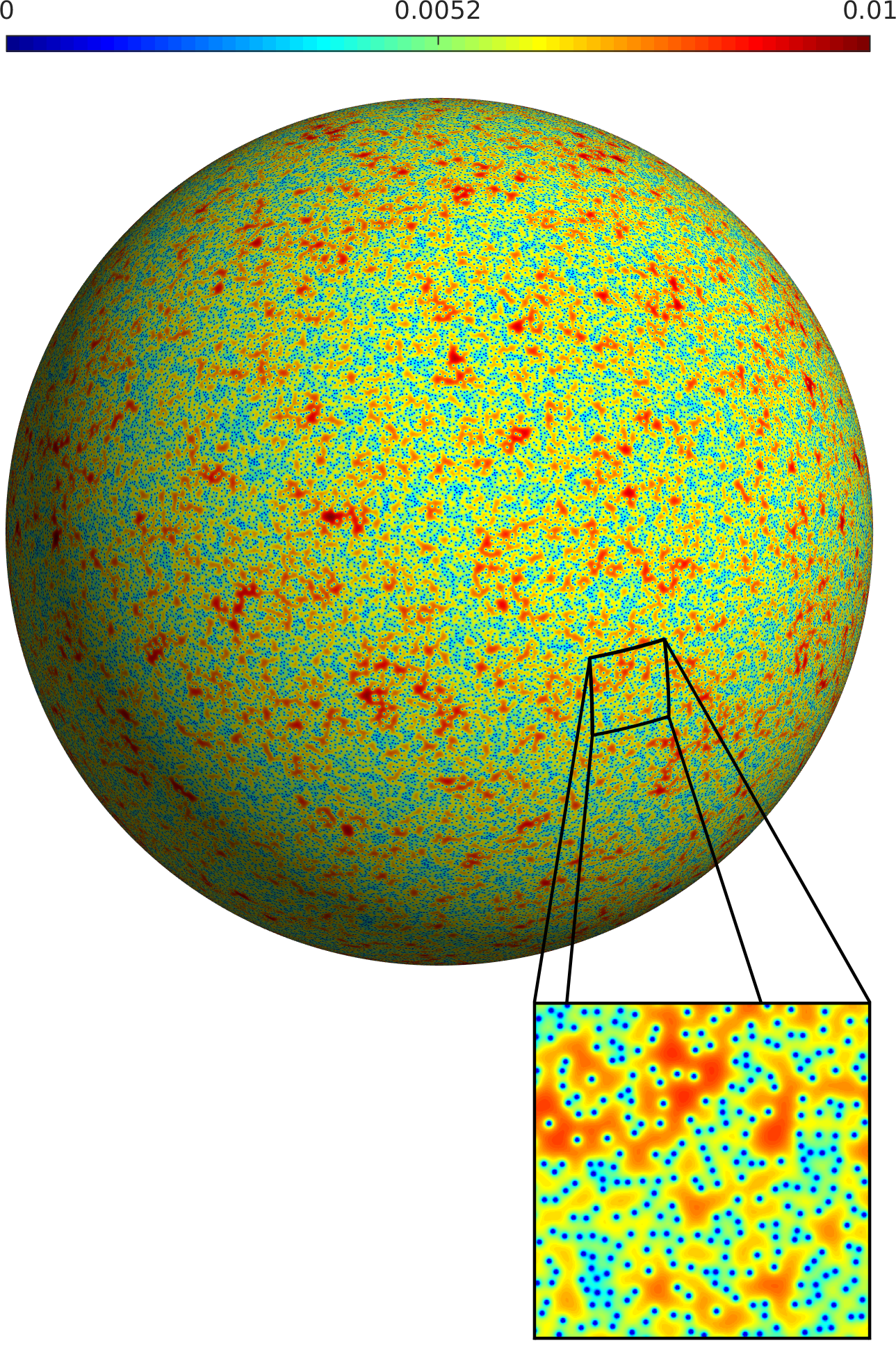}
    \caption{MFPT $\vbar$ plotted just inside the unit sphere for
    an example with $N = 100\, 000$ random well-separated patches of radius $\varepsilon \approx 0.00141$. The
  integral equation associated with this problem was solved in $63$
  minutes on a 60-core workstation, to an $L^2$ residual error of
approximately $2.2 \times 10^{-8}$. Further details are given in Section
\ref{sec:numexamples}.}
    \label{fig:randpts1e5}
\end{figure}

\section{Analytical setup}\label{sec:analyticalsetup}

Our approach to solving the exterior and interior problems \eqref{eq:extBVP}
and \eqref{eq:intBVP} uses a representation of each solution as
an integral involving the corresponding Neumann Green's function. 
This representation leads
to an integral equation, and the scalar quantity of interest - $J$ or
$\mu$ - can be calculated directly from its solution.

\subsection{Neumann Green's functions for the sphere}\label{sec:greens}

Let us first consider the exterior Neumann problem:
\begin{align}\label{eq:extneu}
  \begin{cases}
    \Delta u = 0 & x \in \mathbb{R}^n \backslash \wb{\Omega} \\
    \frac{\partial u}{\partial n} = g & x \in \partial \Omega \\
    u(x) \to 0 & |x| \to \infty.
  \end{cases}
\end{align}
Here $\Omega$ is a bounded domain, and $g$ a given continuous function on
$\partial \Omega$. This
problem has a unique solution, and if $\Omega$ is the unit ball in
$\mathbb{R}^3$, it may be obtained using the {\em exterior Neumann Green's
function} $G_E(x,x')$, which is known analytically \cite{koshlyakov64,kellogg53}.
$G_E$ is symmetric, and
satisfies
\begin{align}\label{eq:Gextproperties}
  \begin{cases}
    -\Delta G_E(x,x') = 4 \pi \delta(x-x') &x,x' \in \mathbb{R}^3 \backslash
    \Omega \\
    \frac{\partial}{\partial n_{x'}} G_E(x,x') = 0 &x \in \mathbb{R}^3
    \backslash \Omega, \, x' \in \partial \Omega, \, x \neq x', \\
  \end{cases}
\end{align}
with $G_E(x,x') = \OO\paren{|x|^{-1}}$ as $|x| \to \infty$ for fixed $x' \in
\mathbb{R}^3 \backslash \Omega$. It can be shown, using Green's
second identity, that
\begin{equation}\label{eq:extneusoln}
  u(x) = \frac{1}{4 \pi} \int_{\partial \Omega} G_E(x,x') g(x') \, dS(x')
\end{equation}
solves the exterior Neumann problem \eqref{eq:extneu}. When $x' \in
\partial \Omega$, $G_E$ is given explicitly by
\begin{equation}\label{eq:Gextoffsurface}
  G_E(x,x') = \frac{2}{|x-x'|} + \log \left( \frac{|x| - x \cdot x'}{1 - x \cdot x' + |x -
  x'|} \right).
\end{equation}
If, in addition, $x \in \partial \Omega$, then
\begin{equation}\label{eq:Gextonsurface}
  G_E(x,x') = \frac{2}{|x-x'|} - \log \left( \frac{2}{|x-x'|} \right) -
\log \left(1 + \frac12 |x-x'| \right).
\end{equation}

The interior Neumann problem is given by
\begin{align}\label{eq:intneu}
  \begin{cases}
    \Delta v = 0 & x \in \Omega  \\
    \frac{\partial v}{\partial n} = g & x \in \partial \Omega,
  \end{cases}
\end{align}
where $\Omega$ is a bounded domain and $g$ is a continuous function defined on the 
boundary, with the additional constraint that $g$
must satisfy the consistency condition
\[\int_{\partial \Omega} g \, dS = 0.\]
This problem has a solution which is unique up to an additive constant.
The consistency condition precludes the existence of an interior Green's function
with zero Neumann data. Rather, for $\Omega$
the unit ball in $\mathbb{R}^3$, we have an {\em interior Neumann Green's
function} $G_I(x,x')$, also known analytically
\cite{koshlyakov64,kellogg53}. It is again
symmetric and satisfies
\begin{align}\label{eq:Gintproperties}
  \begin{cases}
    -\Delta G_I(x,x') = 4 \pi \delta(x-x') &x,x' \in \Omega \\
    \frac{\partial}{\partial n_{x'}} G_I(x,x') = -1 &x \in
    \wb{\Omega}, \, x' \in \partial \Omega, \, x \neq x'. \\
  \end{cases}
\end{align}
As before,
\begin{equation}\label{eq:intneusoln}
  v(x) = \frac{1}{4 \pi} \int_{\partial \Omega} G_I(x,x') g(x') \, dS(x')
\end{equation}
solves the interior Neumann problem \eqref{eq:intneu}. When $x' \in
\partial \Omega$, $G_I$ is given by
\begin{equation}\label{eq:Gintoffsurface}
  G_I(x,x') = \frac{2}{|x-x'|} + \log \left( \frac{2}{1 - x \cdot x' +
  |x - x'|} \right).
\end{equation}
If, in addition, $x \in \partial \Omega$, this reduces to
\begin{equation}\label{eq:Gintonsurface}
  G_I(x,x') = \frac{2}{|x-x'|} + \log \left( \frac{2}{|x-x'|} \right) -
\log \left(1 + \frac12 |x-x'| \right).
\end{equation}
This is the same as \eqref{eq:Gextonsurface} except for the
sign of the second term. In other words,
the restrictions of the interior and exterior Green's functions
to the boundary $\partial \Omega$ are nearly identical.

The following lemma, which we will require in the next section, follows
from the second property in \eqref{eq:Gintproperties} and the symmetry
of $G_I$.

\begin{lemma}\label{lem:intsigma}
  Let $\Gamma$ be an open subset of $\partial \Omega$ and let $\sigma$
  be continuous on $\Gamma$. Then for $x \in \partial \Omega \backslash
  \bar{\Gamma}$,
  \[\frac{\partial}{\partial n_x} \int_{\Gamma} G_I(x,x')
  \sigma(x') \, dS(x') = -\int_{\Gamma} \sigma(x') \, dS(x').\]
\end{lemma}

\begin{figure}[ht]
  \centering
    \hspace*{3cm}\includegraphics[width=0.8\linewidth]{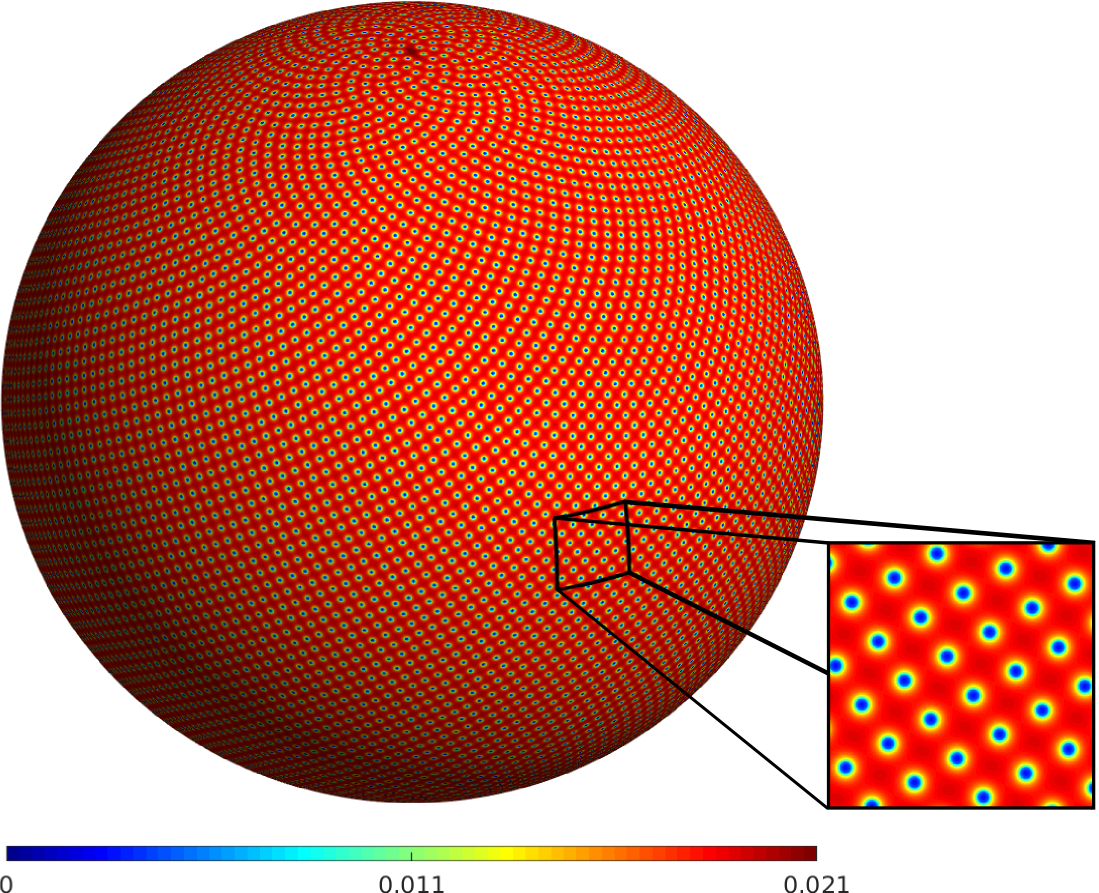}
    \caption{MFPT $\vbar$ plotted just inside the unit sphere for
    an example with $N = 10\, 000$ uniformly distributed patches of radius $\varepsilon
    \approx 0.00447$. The
  integral equation associated with this problem was solved in $114$
  seconds on a 60-core workstation, and in $15$ minutes on a four-core,
  eight-thread laptop, to an $L^2$ residual error of
approximately $6.4 \times 10^{-8}$. Further details are given in Section
\ref{sec:numexamples}.}
    \label{fig:fibopts1e4}
\end{figure}

\subsection{The narrow capture problem}\label{sec:extsetup}

We turn now to the narrow capture problem, which is
the simpler of the two. We first modify the BVP
\eqref{eq:extBVP} by defining $u = 1 - \ubar$, so that solutions decay
as $|x| \to \infty$. 
The function $u$ satisfies the modified equations
\begin{align}\label{eq:extBVP2}
  \begin{cases}
    \Delta u = 0 & x \in \mathbb{R}^3 \backslash \Omega \\
    u = 1 & x \in \Gamma_A \\
    \frac{\partial u}{\partial n} = 0 & x \in \Gamma_R \\
    u(x) \to 0 & |x| \to \infty. \\
  \end{cases}
\end{align} 

Let us denote the unknown Neumann data on $\Gamma_A$ by 
$\sigma(x')$.  Then
 \eqref{eq:extneusoln} implies that for $x \in \mathbb{R}^3
\backslash \wb{\Omega}$, we have
\begin{equation}\label{eq:urep}
  u(x) = \frac{1}{4 \pi} \int_{\Gamma_A} G_E(x,x') \frac{\partial
u}{\partial n}(x') \, dS(x') \equiv \int_{\Gamma_A} G_E(x,x') \sigma(x')
\, dS(x').
\end{equation}
By analogy with classical potential theory, we refer to this as a {\em
single layer potential} representation with density $\sigma$ supported
on $\Gamma_A$. Since the dominant singularity of the kernel $G_E$ is that of
the free-space Green's function for the Laplace equation, this single layer potential is continuous up to $\partial \Omega$.
Taking the limit as $x \rightarrow \Gamma_A$ and using
the second condition in \eqref{eq:extBVP2}, we obtain the first-kind
integral equation
\begin{equation}\label{eq:extinteqn}
  \int_{\Gamma_A} G_E(x,x') \sigma(x') \, dS(x') = f(x), \quad x \in \Gamma_A,
\end{equation}
where $f(x) \equiv 1$, with the weakly singular kernel $G_E$. Assuming that we can solve
\eqref{eq:extinteqn} for $\sigma$, it follows that
$u(x)$, given by \eqref{eq:urep}, is the solution to 
\eqref{eq:extBVP2}, and that $\ubar = 1-u$ solves
\eqref{eq:extBVP}. Furthermore, since $\sigma \equiv \frac{\partial
u}{\partial n} \equiv -\frac{\partial \ubar}{\partial n}$ on
$\Gamma_A$, the total flux $J$ from \eqref{eq:flux} will be given by
\[J = -I_\sigma\]
where we have introduced the shorthand 
\begin{equation}
\label{Isigmadef}
I_\sigma := \int_{\Gamma_A} \sigma \, dS. 
\end{equation}

We will not prove the existence of a solution to \eqref{eq:extinteqn},
but sketch a possible approach.
If we replace the kernel $G_E$ in \eqref{eq:extinteqn} with its first term
$\frac{2}{|x-x'|}$, which is the free-space Green's function for the Laplace
equation (up to a constant scaling factor), we obtain the
first-kind integral equation for the Dirichlet problem 
on an open surface, which we can denote in operator form by 
\[ \calS_0 \sigma = f. \] 
This is
a well-studied problem, which has a unique solution in the Sobolev space
$H^{-\frac12}(\Gamma_A)$ given data in 
$H^{\frac12}(\Gamma_A)$ \cite{Stephan1987}. 
Writing the full single layer potential operator in the form
$\calS_0 + K$, where $K$ is a 
compact pseudodifferential operator of order $-2$,
we may rewrite 
\eqref{eq:extinteqn} in the form of a Fredholm integral equation of the second kind:
\begin{equation}
 (I + \calS_0^{-1}K) \sigma = \calS_0^{-1} \, f. 
\label{preconsinglepatch}
\end{equation}
Thus, to prove existence and uniqueness for the single patch equation,
one can apply the Fredholm alternative to 
\eqref{preconsinglepatch}.
That is, one need only show that the homogenous 
version of the single patch equation has no nontrivial solutions.
This is straightforward to prove when
$\varepsilon$ is sufficiently small, since the norm of $K$ goes to zero
as $\varepsilon$ goes to zero and the corresponding Neumann series 
converges.
We conjecture that the result holds for any $\varepsilon$.

\subsection{The narrow escape problem}\label{sec:intsetup}

The analytical formulation of the narrow escape problem is somewhat more complicated
than that of the narrow capture problem, largely because of the
non-uniqueness of the interior Neumann problem, but it leads to a similar integral equation. 
We first recast the Poisson problem
\eqref{eq:intBVP} as a Laplace problem
with inhomogeneous boundary conditions. 
Assume that $v$ satisfies 
\begin{align}\label{eq:intBVP2}
  \begin{cases}
    \Delta v = 0 & x \in \Omega \\
    v = 1 & x \in \Gamma_A \\
    \frac{\partial v}{\partial n} = D & x \in \Gamma_R,
  \end{cases}
\end{align}
for some non-zero constant $D$. Then $\vbar$ given by
\begin{equation}\label{eq:vtilde}
  \vbar = \frac{v - 1}{3D} + \frac{1 - |x|^2}{6}
\end{equation}
solves \eqref{eq:intBVP}. We will therefore seek a method to produce a
solution of \eqref{eq:intBVP2} for some $D \neq 0$.

\begin{figure}[ht]
  \centering
    \hspace*{3cm}\includegraphics[width=0.8\linewidth]{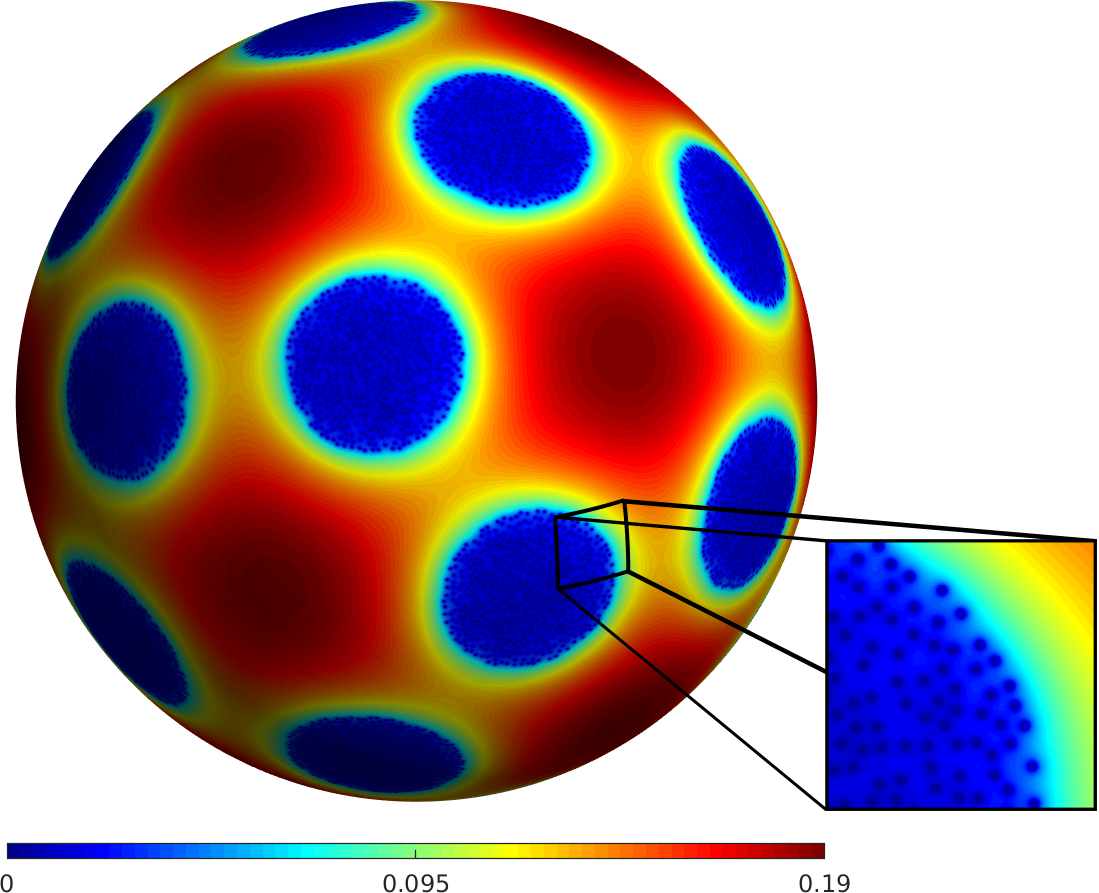}
    \caption{MFPT $\vbar$ plotted just inside the unit sphere for
    an example with $N = 10\, 000$ random, clustered patches of radius $\varepsilon
    \approx 0.0035$. The
  integral equation associated with this problem was solved in $269$
  seconds on a 60-core workstation, and in $35$ minutes on a four-core,
  eight-thread laptop, to an $L^2$ residual error of
approximately $6.5 \times 10^{-8}$. Further details are given in Section
\ref{sec:numexamples}.}
    \label{fig:cluspts1e4}
\end{figure}

\begin{lemma}
Let
\begin{equation}\label{eq:intBVPansatz}
  v(x) = \int_{\Gamma_A} G_I(x,x') \sigma(x') \, dS(x'),
\end{equation}
where $\sigma$ satisfies the first-kind integral
equation
\begin{equation}\label{eq:intinteqn}
  \int_{\Gamma_A} G_I(x,x') \sigma(x') \, dS(x') = 1
\end{equation}
for $x \in \Gamma_A$. Then $v$ solves \eqref{eq:intBVP2} with $D = -I_\sigma$, for $I_\sigma$
defined as in \eqref{Isigmadef}, and $I_\sigma \neq 0$.
\end{lemma}

\noindent {\em Proof}:
The function $v(x)$ is 
harmonic in $\Omega$, and by Lemma \ref{lem:intsigma}, it
satisfies the third condition of \eqref{eq:intBVP2} with $D \equiv
-I_\sigma$, as long as $I_\sigma \neq 0$. Taking $x$ to $\Gamma_A$ and using the continuity of
the single layer potential up to
$\Gamma_A$, we find that $v$ will satisfy the second condition of
\eqref{eq:intBVP2} as long as $\sigma$ satisfies \eqref{eq:intinteqn}.

It remains only to show that if $\sigma$ satisfies 
\eqref{eq:intinteqn}, then
$I_\sigma \neq 0$.
If not, then $v$ given by \eqref{eq:intBVPansatz} satisfies
\eqref{eq:intBVP2} with $D = 0$, as does the constant function $1$. 
It follows from Green's identity
that solutions to \eqref{eq:intBVP2} with the same value of $D$ are unique, so we
must have $v \equiv 1$. The formula
\eqref{eq:Gintoffsurface} for $G_I$ shows that if $|x'| = 1$, then $G_I(0,x')
= 2$, so if $v \equiv 1$ we have
\[ 1 = v(0) = 2 \int_{\Gamma_A} \sigma(x') \, dS(x') = 2 I_\sigma, \]
a contradiction. \hfil \qed \\
The question of the existence of a solution to \eqref{eq:intinteqn} is
analogous to that for \eqref{eq:extinteqn}, which was discussed in
Section \ref{sec:extsetup}.

To calculate the average MFPT $\mu$ directly from
$\sigma$, we plug \eqref{eq:vtilde} into \eqref{eq:avgmfpt2} to obtain
\begin{equation}\label{eq:mufromsigma1}
  \mu = \frac{1}{3 D
  |\partial \Omega|} \int_{\partial \Omega} v \, dS - \frac{1}{3D} +
\frac{1}{15}.
\end{equation}
To calculate $\frac{1}{|\partial \Omega|} \int_{\partial \Omega} v \, dS$, we use the
representation \eqref{eq:intBVPansatz}:
\begin{align*}
  \frac{1}{|\partial \Omega|} \int_{\partial \Omega} v \, dS &= \frac{1}{|\partial \Omega|} \int_{\partial \Omega} \int_{\Gamma_A}
G_I(x,x') \sigma(x') \, dS(x') \, dS(x) \\
&= \int_{\Gamma_A} \sigma(x') \paren{ \frac{1}{|\partial \Omega|} \int_{\partial \Omega}
G_I(x,x') \, dS(x) } \, dS(x').
\end{align*}
A calculation using the explicit form \eqref{eq:Gintonsurface} of $G_I$ gives
\[\frac{1}{|\partial \Omega|} \int_{\partial \Omega} G_I(x,x') \, dS(x) = 2\]
for any $x' \in \partial \Omega$.
We therefore have
\[\frac{1}{|\partial \Omega|} \int_{\partial \Omega} v \, dS = 2 I_\sigma.\]
Plugging this into \eqref{eq:mufromsigma1} and replacing $D$ by
$-I_\sigma$ gives
\begin{equation}\label{eq:mufromsigma2}
  \mu = \frac{1}{3 I_\sigma} - \frac35.
\end{equation}

\section{A multiple scattering formalism}\label{sec:manypatchsetuptools}

We have shown that the solutions of the two boundary value problems of
interest, as well the associated scalars $J$ and $\mu$, may be
obtained by solving \eqref{eq:extinteqn} and \eqref{eq:intinteqn},
respectively,
on the collection of absorbing patches. These integral equations differ only by the sign of one 
term in their respective kernels, as seen in Section \ref{sec:greens}.
Since our treatment of the two cases is the same,
we drop the subscripts on $G_E$ and $G_I$, and discuss the solution of
\[ \int_{\Gamma_A} G(x,x') \sigma(x') \, dS(x') = 1 \quad x \in \Gamma_A,\] 
where $\sigma$ is an unknown density on $\Gamma_A$. Letting
$\Gamma_A = \cup_{i=1}^N \Gamma_{i}$, where $\Gamma_{i}$ is the $i$th
patch, and letting $\sigma_i$ be the restriction of $\sigma$ 
to $\Gamma_{i}$, we write this equation in the form
\begin{equation}\label{eq:inteqns}
  \sum_{j=1}^N \int_{\Gamma_{j}} G(x,x') \sigma_j(x') \, dS(x') = 1 \quad x \in
\Gamma_{i}, \, i = 1,\ldots,N.
\end{equation} 

For the sake of simplicity, we assume that each patch has the same radius 
$\varepsilon$. We also assume that the patches are {\em well-separated},
in the sense that the distance between the centers of any two patches
in arc length along the surface of the sphere is at least $3 \varepsilon$.
That is, any two patches are separated by a distance greater than or
equal to their own radius. 
For $x \in \Gamma_i$, we define $\calS_{ij}$ by
\[(\calS_{ij} \sigma_j)(x) := \int_{\Gamma_j} G(x,x') \sigma_j(x')
\, dS(x'). \]
More specifically, we define each such operator in a coordinate system
fixed about the center of $\Gamma_j$. Since all the patches have the same radius, the operators $\calS_{ii}$
are therefore identical, and we denote $\calS_{ii}$ by $\calS$. Thus we may rewrite the many-patch integral equation \eqref{eq:inteqns} in the form
\begin{equation}\label{eq:opinteqns}
  \calS \sigma_i + \sum_{j\neq i}^N \calS_{ij} \sigma_j = 1 \quad i = 1,\ldots,N.
\end{equation}
The aim of this section is to reformulate \eqref{eq:opinteqns} as a
Fredholm system of the second kind in an efficient basis.

\begin{definition} \label{def:1pinteqn}
Let $f$ be a smooth function on some patch $\Gamma_i$. 
The \emph{one-patch integral
  equation with data $f$} is defined by
  \begin{equation}\label{eq:onepatinteq}
    \calS \sigma_i = f,
  \end{equation}
  where $\sigma_i$ is an unknown density on $\Gamma_i$.
\end{definition}

\begin{remark} \label{rem:rightside}
  Writing \eqref{eq:opinteqns} in the form
\[
\calS \sigma_i  =  1 - \sum_{j\neq i}^N \calS_{ij} \sigma_j, 
\]
and observing that $\calS_{ij} \sigma_j$ is a smooth function for
$\Gamma_j$ well-separated from $\Gamma_i$, we see that each $\sigma_i$
satisfies a one-patch integral equation with smooth data. 
Conversely, if $\sigma_1,\ldots,\sigma_N$ satisfy \eqref{eq:opinteqns}, then
each $\calS \sigma_i$ is smooth on $\Gamma_i$.
\end{remark}

It is convenient to
make use of an orthonormal basis $\{ q_1,q_2,\dots \}$ of smooth functions
on each patch, so that for smooth $f$ on $\Gamma_i$ we have
\begin{equation} \label{fsynth}
  f(x) = \sum_{n=1}^\infty \hat{f}_n q_n(x),
\end{equation}
in the usual $L^2$ sense, with 
\[ \hat{f}_n = \int_{\Gamma_i} f(x) q_n(x) \, dx.
\]
We postpone until Section \ref{sec:zernike} a discussion of our
particular choice of the basis $\{q_n\}$, which will be constructed 
using Zernike polynomials.
We will denoted by $\hat{f}^K$ the 
vector of the first $K$ coefficients:
\[ \hat{f}^K = (\hat{f}_1,\hat{f}_2,\ldots,\hat{f}_K)^T.
\] 

\begin{definition}  \label{projectors}
Let $f$ be a smooth function on $\Gamma$ defined by 
\eqref{fsynth}, with $\hat{f}$, $\hat{f}^K$ computed as above. 
The {\em projection}
operators $\proj$ and $\projK$ are defined by
\[ \paren{\proj[f]}_n = \hat{f}_n,  \]
with $\projK$ defined in the same manner for $n \leq {K}$.
The {\em synthesis} operators $\syn$ and $\synK$ are
defined by
\[ \syn[\hat{f}](x) = \sum_{n=1}^\infty \hat{f}_n q_n(x),\quad
  \synK[\hat{f}_{K}](x) = \sum_{n=1}^{K} f_n q_n(x).
\]
\end{definition}
$\proj$ and $\projK$ are left inverses of 
$\syn$ and $\synK$, respectively.

Finally, we define $b_n$ to be the solution of 
the one-patch integral equation
with data given by the basis element $q_n$:
\begin{equation}\label{eq:onepatchinteqn}
  b_n = \calS^{-1} q_n.
\end{equation}
Thus, if a smooth function $f$ on $\Gamma_i$ is expanded
as $f = \sum_{n=1}^\infty \hat{f}_n q_n$, then the
solution of the one-patch integral equation with data $f$ is given by
$\calS^{-1} f =  \sum_{n=1}^\infty \hat{f}_n b_n$. This
motivates the following definition.

\begin{definition}\label{def:Bop}
  We denote the solution operator of the one-patch integral equation in the
  basis $\{q_n\}$ by
  \[\calB = \calS^{-1} \syn.\]
  For $\hat{f} = \{\hat{f}_1,\hat{f}_2,\ldots\}$ and $f(x) =
  \sum_{n=1}^\infty \hat{f}_n q_n(x)$, $\calB$ satisfies
  \[\calB[\hat{f}](x) = \sum_{n=1}^\infty \hat{f}_n b_n(x).\]
  We denote the solution operator of the one-patch integral equation in the
  truncated basis $\{q_n\}_{n=1}^K$ by
  \[\calBK = \calS^{-1} \synK.\]
  For $\hat{f} = (\hat{f}_1,\hat{f}_2,\ldots\,\hat{f}_K)$ and $f(x) =
  \sum_{n=1}^K \hat{f}_n q_n(x)$, $\calB_K$ satisfies
  \[\calBK[\hat{f}](x) = \sum_{n=1}^K \hat{f}_n b_n(x).\]
\end{definition}

Note that the construction of $\calB$ requires
solving the one-patch integral equations with data $q_1,q_2,\ldots$ to
obtain $b_1,b_2,\ldots$, and that the 
construction of $\calBK$ requires solving the first $K$ of these
equations. For a fixed patch radius
$\varepsilon$, these solutions are universal and do not 
depend on the number or arrangement of patches in the full problem.

Given $\calB$, we are now able to rewrite the integral equation 
\eqref{eq:opinteqns} as a well-conditioned Fredholm system of the second
kind in the basis $\{q_n\}$. On $\Gamma_i$, we {\em define} 
a function $f_i$ by
\[ f_i = \calS \sigma_i. \]
Substituting into \eqref{eq:opinteqns}, we have 
\[f_i + \sum_{j\neq i}^N \calS_{ij} \calS^{-1} f_j = 1
\quad i = 1,\ldots,N. \]
To transform to the basis $\{q_n\}$, we write
$f_i$ in the form $f_i = \syn \hat{f}_i$
and multiply on the left by $\proj$ to obtain
\begin{equation}\label{eq:opinteqnslim*}
  \hat{f}_i + \proj \sum_{j\neq i}^N \calS_{ij}\calB \hat{f}_j = 
  \proj \, 1 \quad i = 1,\ldots,N.
\end{equation}
Since the patches $\Gamma_i$ and $\Gamma_j$ are well-separated, 
$\proj \calS_{ij} \calB$ is a compact operator for $i \neq j$, 
so that \eqref{eq:opinteqnslim*} is a Fredholm system of the second kind.
The corresponding truncated system takes the form 
\begin{equation}\label{eq:opinteqns*}
\hat{f}_i^K + \projK \sum_{j\neq i}^N 
\calS_{ij}\calBK \hat{f}_j^K = 
\projK \, 1 \quad i = 1,\ldots,N,
\end{equation}
where we have used the approximation $f_i \approx \synK \hat{f}_i^K$.

\begin{remark}\label{rem:advantages}
We refer to the approach described above as a {\em multiple scattering}
formalism by analogy with the problem of wave scattering from multiple
particles in a homogeneous medium. In the language of scattering theory,
one would say that for the $i$th patch, the boundary
data is the known data ($\calS \sigma_i = 1$), perturbed by 
the potential ``scattered" from all other patches, namely
$\sum_{j\neq i}^N \calS_{ij} \sigma_j$.
Solving the system \eqref{eq:opinteqns} corresponds to determining how the 
collection of uncoupled single patch solutions
$\calS \sigma_i = 1$ needs to be perturbed to account for the 
``multiple scattering" effects. 

The approach developed above, where 
$f_i = \calS \sigma_i$ are the unknowns, has many advantages
  over solving \eqref{eq:opinteqns} directly, even with ${\calS}^{-1}$ as 
  a left preconditioner. By working in the spectral basis, we avoid the
  need to discretize $\sigma_i$ on each patch, the number
  of degrees of freedom per patch is significantly reduced, and the 
  linear system is a well-conditioned Fredholm equation of the second kind.
\end{remark}

\begin{remark}\label{rem:sigrep}
  The original unknowns $\sigma_i$ may be recovered from the
  solution of \eqref{eq:opinteqnslim*} or \eqref{eq:opinteqns*} 
  using the formula
  \begin{equation}\label{eq:recoversig}
    \sigma_i = \calB \hat{f}_i \approx \calB^K \hat{f}_i^K.
  \end{equation}
  Thus, we may think of the unknowns
  $\hat{f}_i$ as a representation of the unknown density $\sigma_i$ in the 
  basis $\{b_n\}$.
\end{remark}

We turn now to 
the construction of an orthonormal basis $\{ q_n \}$ for 
smooth functions on a patch, the  
construction of the singular solutions $b_n = \calS^{-1} q_n$, and
the efficient solution of the discretized multiple scattering system 
\eqref{eq:opinteqns*}.

\section{A basis for smooth functions on a patch} \label{sec:zernike}

It is well-known that the Zernike polynomials are a
spectrally accurate, orthogonal basis for smooth functions on the
disk.
For a thorough discussion of these functions, we refer the reader
to \cite{boyd11}. Here, we simply summarize their relevant properties.

The Zernike polynomials on the unit disk $0 \leq r \leq 1$, $0 \leq
\theta < 2 \pi$ are given by
\begin{align*}
  \begin{cases}
    Z_n^m(r,\theta) &= R_n^m(r) \cos(m \theta) \\
    Z_n^{-m}(r,\theta) &= R_n^m(r) \sin(m \theta),
  \end{cases}
\end{align*}
with $0 \leq m < \infty$, $m \leq n < \infty$, and 
\[R_n^m(r) = (-1)^{(n-m)/2} r^m P_{(n-m)/2}^{m,0}(1 - 2 r^2), \]
where $P_n^{\alpha,\beta}(x)$ is a Jacobi polynomial on
$[-1,1]$. The Jacobi polynomials are orthogonal on $[-1,1]$ with respect to the
weight function $(1-x)^\alpha (1+x)^\beta$. Thus, for fixed
$m$, the functions $R_n^m(r)$ are orthogonal on $[0,1]$ with respect to
the weight function $r$. This gives the orthogonality relation
\begin{equation} \label{eq:zorthog}
  \int_0^{2 \pi} \int_0^1 Z_{n_1}^{m_1}(r,\theta)
Z_{n_2}^{m_2}(r,\theta) r \, dr \, d\theta = \frac{(1 + \delta_{m_1,0}) \pi}{2
n_1 + 2} \delta_{n_1,n_2} \delta_{m_1,m_2}.
\end{equation}

The natural truncation of this basis is to fix a cutoff mode $M$ in both
the radial and angular variables, and to let $0 \leq m \leq n \leq M$. 
This yields $K = (M+1)(M+2)/2$ basis functions. 
To use this basis on a generic patch $\Gamma_i$, 
we define a polar coordinate system
$(r,\theta)$ about the patch center, for which $r$ is the distance in
arc length along the sphere from the center, and $\theta$ is the
polar angle. We rescale the radial variable 
from $[0,1]$ to $[0,\varepsilon]$, transforming the Zernike
polynomials to functions on $\Gamma_i$. Finally, the basis functions
$q_1,\ldots,q_K$ discussed in Section \ref{sec:manypatchsetuptools} can be defined as the scaled Zernike
polynomials up to mode $M$.

>From the orthogonality relation \eqref{eq:zorthog},
the projection operators $\proj$ and $\projK$ 
are obtained as
normalized inner products against Zernike polynomials in polar coordinates.
This {\em Zernike transform} can be implemented
numerically using
a tensor product quadrature with a Gauss-Legendre rule in the radial
variable and a trapezoidal rule in the angular variable. 
The number of grid points required to obtain the exact Zernike
coefficients of a function in the space spanned by 
$q_1,\dots,q_K$ is $\OO(K)$; we denote this number by $K^*$. We refer to these points as
the {\em Zernike sampling nodes} $x_1^z,\ldots,x_{K^*}^z$
(see \cite{boyd11} for further details).

\begin{remark}\label{rem:zerniketrunc}
Rewriting \eqref{eq:opinteqns*} in the form
  \begin{equation}\label{eq:opinteqns2*}
    \hat{f}_i^K = 
  \projK \, \left(1 - \sum_{j\neq i} \calS_{ij}\calBK \hat{f}_j^K \right),
  \end{equation}
  we see that the truncation error compared with \eqref{eq:opinteqnslim*} depends on how well 
the smooth function
\[1 - \sum_{j\neq i} \calS_{ij}\calBK \hat{f}_j^K\]
is represented in the space spanned by
$q_1,\ldots,q_K$. 
In the one-patch case, the summation term vanishes, and $K=1$ is sufficient. 
For multiple patches, the choice of $K$ depends largely on 
how well-separated the patches are.
Since the Zernike basis is spectrally accurate, $M$ grows only logarithmically 
with the desired precision. In practice, {\em a posteriori} estimates 
are easily obtained for any fixed configuration by inspection of the 
decay of the Zernike coefficients $\hat{f}_i^K$ in the computed solution.
\end{remark}

\section{Informal description of the one-patch solver}
\label{sec:onepatchprelim}

While the details of our solver for the one-patch integral equation
\[ \calS \sigma_i = f \]
are deferred to Section \ref{sec:onepatch}, we outline the general 
approach here.
First, we note that in the absence of curvature 
(i.e. a flat disk on a half-space) and with the associated terms of the
Green's function removed, the solution $\sigma_i$ 
is known to have a square root singularity at the disk edge
\cite{Bernoff2018,Sneddon1966,Fabrikant1989,Stephan1987,costabelsing}.
In our case, we will explicitly include this square root singularity in
the representation of $\sigma_i$, but also allow for weaker
singularities - which
we have observed and will demonstrate in Section
\ref{sec:sigsingularity} - by using a
discretization that is adaptively refined toward the edge 
$\partial \Gamma_i$.

Assume then that we have discretized 
the patch $\Gamma_i$ using a suitable polar mesh with 
$n_f$ {\em fine grid points}, denoted by 
$x_{i,1}^f,\ldots,x_{i,n_f}^f$. The fine grid points for different
patches are identical relative to the coordinate systems of their own
patches. We denote the corresponding samples
of the right-hand side $f$ and $\sigma_i$ by
\[
\begin{aligned}
  \vec{f} &= (f(x_{i,1}^f),\ldots,f(x_{i,n_f}^f))^T, \\
 \vec{\sigma}_i &= ((\vec{\sigma}_i)_1,\ldots,(\vec{\sigma}_i)_{n_f})^T 
\approx
(\sigma_i(x_{i,1}^f),\ldots,\sigma_i(x_{i,n_f}^f))^T. 
\end{aligned} 
\]
We assume that $\calS$ is discretized to high-order accuracy by a matrix
$S$ with
\begin{equation}\label{eq:finegrdquadraturea}
  \calS[\sigma_i](x_{i,k}^f) 
  \approx \sum_{l=1}^{n_f} S(k,l) (\vec{\sigma_i})_l,
\end{equation}
so that the discretized system takes the form
\begin{equation}\label{eq:onepatchdiscreteeqn}
  S \vec{\sigma}_i = \vec{f}.
\end{equation}

We will also require a set of quadrature weights, denoted by
$w_1^f,\ldots,w_{n_f}^f$ and identical for each patch, that permit the accurate integration
over $\Gamma_i$ of the product of an arbitrary smooth function with the 
discretized density $\vec{\sigma}_i$, taking into account the fact that
$\sigma_i$ has an edge singularity.
That is, we assume that
\begin{equation}\label{eq:finegrdquadrature}
  \int_{\Gamma_i} g(x) \sigma_i(x) \, dS(x) \approx \sum_{l=1}^{n_f}
  g(x_l^f) (\vec{\sigma}_i)_l w_l^f
\end{equation}
for any smooth $g$, with high-order accuracy. 
In the next section, we
will use this quadrature to discretize the operators $\calS_{ij}$.

The solutions of the $K$ one-patch integral
equations \eqref{eq:onepatchinteqn} may be obtained in a precomputation,
after which we have access
to the functions $b_1,\ldots,b_K$ sampled on the fine grid. 
We assemble these
functions into an $n_f \times K$ matrix $B$ with
\[B(n,m) = b_m(x_n^f).\]
$B$ is then the discretization of the operator $\calBK$, mapping the
first $K$ Zernike coefficients of a smooth function to the solution of the
corresponding one-patch integral equation sampled on the fine grid. 
If we denote by $Q$ the discretization of the synthesis operator 
$\synK$ as an $n_f
\times K$ matrix,
\[Q(i,j) = q_j(x_i^f),\]
then we have, as in Definition \ref{def:Bop},
\[SB = Q.\]
In short, the precomputation amounts to solving this matrix system for
$B$.

\section{Discretization of the multiple scattering system}\label{sec:manypatchsetup}

We return now to the multiple scattering system \eqref{eq:opinteqns*}.
The unknowns on $\Gamma_i$ are defined in the truncated Zernike basis as
$\hat{f}_i^K$. We will need as intermediate variables the fine grid
samples of $\sigma_i(x)$. From Remark \ref{rem:sigrep}, we define the
sampling vector $\vec{\sigma_i}$ by
\[ \vec{\sigma_i} = B \hat{f}_i^K \approx \calB^K \hat{f}_i^K. \] 
In order to discretize the integral operators $\calS_{ij}$ 
for $i \neq j$, we note that $G(x,x')$ is smooth for $x \in \Gamma_i$,
$x' \in \Gamma_j$, and
use the quadrature \eqref{eq:finegrdquadrature}. This yields
\begin{equation}\label{eq:Sijquad}
  \int_{\Gamma_{j}} G(x,x') \sigma_j(x') \, dS(x') \approx
\sum_{l=1}^{n_f} G(x,x_{j,l}^f) (\vec{\sigma_j})_l w_l^f.
\end{equation}
Setting $x = x_{i,k}^z$ to be the $k$th Zernike sampling node on
$\Gamma_{i}$, we define the matrix $S_{ij}$ by 
\[S_{ij}(k,l) = G(x_{i,k}^z,x_{j,l}^f) w_l^f. \] 
Thus, $S_{ij}$ maps a density sampled on the fine grid on $\Gamma_j$ to the
smooth field it induces at the Zernike sampling nodes on $\Gamma_i$.
Lastly, we discretize the truncated Zernike transform $\projK$ as a $K
\times K^*$ matrix $P$ using the trapezoidal-Legendre scheme described in
Section \ref{sec:zernike}. 

\begin{definition}
The {\em discrete Zernike transform} $P$ is defined to be the
mapping of a smooth function
sampled on the $K^*$ Zernike sampling nodes to its $K$ Zernike coefficients.
\end{definition}

We can now write the multiple scattering system
\eqref{eq:opinteqns*} in a fully discrete form,
\begin{equation}
\label{eq:opinteqnsdisc}
\hat{f}_i^K + P  \sum_{j \neq i} S_{ij} B 
\hat{f}_j^K   = P \vec{1} \quad i=1,\ldots,N,
\end{equation}
where $\vec{1}$ is the vector of length $K^*$ with all entries
equal to $1$. Since $P \in \RR^{K \times K^*}$,
$S_{ij} \in \RR^{K^* \times n_f}$, and
$B \in \RR^{n_f \times K}$, this is
a linear system of dimensions $K N \times K N$, with
$K << n_f$ degrees of freedom per patch. As a discretization of a
Fredholm system of the second kind, it is amenable to rapid solution
using an iterative method such as GMRES \cite{GMRES}.

We now describe how to calculate the constants $J$ and $\mu$ from the
solution of \eqref{eq:opinteqnsdisc}.
We saw in Sections \ref{sec:extsetup} and
\ref{sec:intsetup} that these can be computed directly from $I_\sigma =
\sum_{i=1}^N \int_{\Gamma_{i}} \sigma_i \, dS$. 
Using the fine grid quadrature \eqref{eq:finegrdquadrature}, we have
\begin{equation}\label{eq:getI}
  I_\sigma = \sum_{i=1}^N
\int_{\Gamma_{i}} \sigma_i \, dS \approx \sum_{i=1}^N \sum_{k=1}^{n_f} (B
\hat{f}_i^K)_k
w_k^f =  (w_1^f,\ldots,w_{n_f}^f) B  \sum_{i=1}^N \hat{f}_i^K.
\end{equation}
Since we may precompute the row vector
$I := (w_1^f,\ldots,w_{n_f}^f) B$ of length $K$,
the cost to compute $I_\sigma$ is $\OO(NK)$.

When the system \eqref{eq:opinteqnsdisc} is solved iteratively,
each matrix-vector product is dominated by the computation
of the ``multiple scattering events''
\begin{equation}\label{eq:mssums}
P \sum_{j \neq i} S_{ij} B \hat{f}_j^K
\end{equation}
for $i=1,\ldots,N$. That is, for each patch $\Gamma_i$, we must compute the
Zernike coefficients of the field induced on that patch by the densities
on all other patches. Note that if we were to calculate the above sums
by simple matrix-vector products, the cost would be
$\OO(n_f K N^2)$. We turn now to the description of a scheme that permits the
computation of these sums using $\OO(K N \log N)$ operations, 
with a constant which depends only on 
the desired precision, but not on $n_f$. 

\section{Efficient representation of outgoing
and incoming fields}\label{sec:outgoingrep}

\begin{figure}[ht]
  \centering
    \includegraphics[width=.4\linewidth]{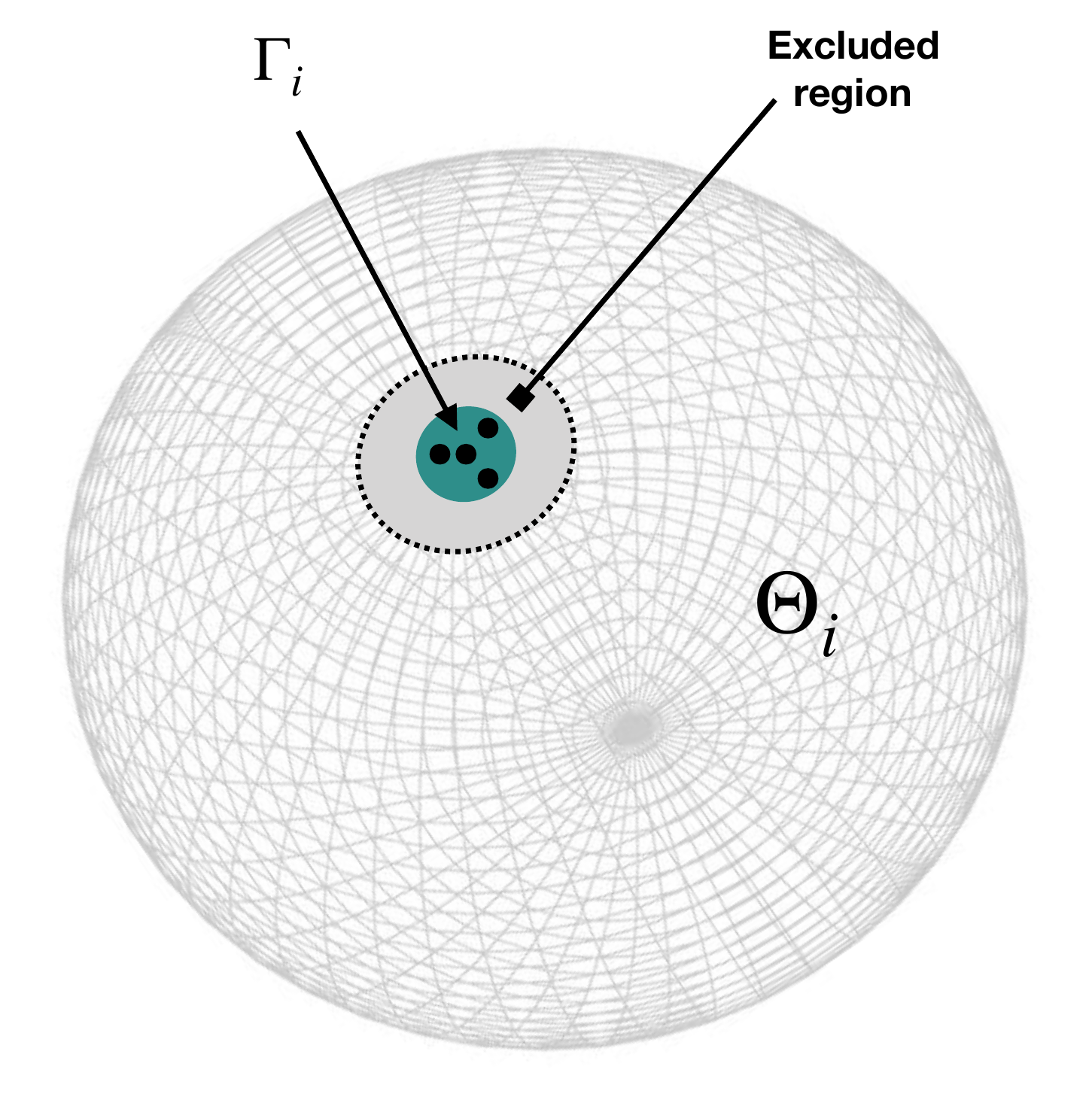}
    \caption{
    For a patch $\Gamma_i$, the far field region $\Theta_i$ is defined as
    the complement on the surface of the sphere of a disk of radius $2
    \varepsilon$, measured in arclength, about the center of $\Gamma_i$.  The black dots in the
    figure represent the subset of the fine grid points used to
    efficiently represent the outgoing field induced by the density
    $\sigma_i$.}
    \label{fig:skelpts}
\end{figure}

Our fast algorithm relies on what is variously referred to as a
compressed, skeletonized, or sparsified representation of the far field
induced by a source density $\sigma_i$ on a single patch $\Gamma_i$ (Fig.
\ref{fig:skelpts}).
We define the far field region $\Theta_i$ for a patch $\Gamma_i$ 
to be the set of points whose distance from the center of $\Gamma_i$ (measured
in arc length along the surface of the sphere) is greater than
$2\varepsilon$. 
In light of our restriction on the minimum patch
separation distance, this ensures that the far field region of a
particular patch contains every other patch.

We start from \eqref{eq:Sijquad}, which was used to define the matrix
$S_{ij}$.  We will show that there is a subset of $p$ fine grid points
with $p << n_f$ and modified source strengths $\vec{\rho}_i =
(\rho_{i,1},\rho_{i,2},\dots,\rho_{i,p})^T$ so that 
\begin{equation}\label{eq:outgoingcompressed}
 \int_{\Gamma_i} G(x,x') \sigma_i(x') \, dS(x') \approx 
  \sum_{l=1}^{n_f} G(x,x_{i,l}^f) (\vec{\sigma}_i)_l w_l^f 
\approx \sum_{m=1}^p
  G(x,x_{i,\pi(m)}^f) \rho_{i,m},
\end{equation}
for any $x \in \Theta_i$. Moreover, there is a stable algorithm
for obtaining this {\em compressed} or
{\em skeletonized outgoing representation}. Here, 
$\pi(m)$ is an indexing function which maps
$\{1,\ldots,p\} \to \{1,\ldots,n_f\}$, and identifies which of the 
original fine grid points are used in the representation. 
The number $p$ represents the numerical rank, to a specified precision,
of the $n_f$ functions $\{G(x,x_{i,l}^f)\}$ on $\Theta_i$. 

\begin{remark}
The existence of such low-rank factorizations is
discussed in detail in \cite{cheng05,Goreinov1997,gu96}.
For the purposes of computation, we will use the
interpolative decomposition (ID) \cite{liberty07,cheng05,id2011}, described
briefly below. 
The ID and related compression schemes are essential and widely used in 
hierarchical, fast algorithms
for applying and inverting dense matrices (see for example
\cite{siva,borm2003hierarchical,corona2014n,fong2009,
gillman2012,gimbutas2002,ho2011,martinsson2005fast,minden2017,ying2004}
and the references therein).
\end{remark}

\subsection{The interpolative decomposition}\label{sec:ID}

We consider a generic patch $\Gamma_i$ and, for simplicity, drop 
the patch index $i$ on all quantities.
We first discretize $\Theta$ on a
training grid $x_1^t,\ldots,x_{n_t}^t$ of $n_t$ points chosen to
be sufficiently fine to
accurately represent smooth functions on $\Theta$. We can then obtain a
matrix $A$ of size $n_t \times n_f$, with entries $A_{jl} =
G(x_j^t,x_l^f)$, so that the $l$th column of $A$ is a discretization of
the function $G(x,x_l^f)$ on the training grid. 
Given a user-specified tolerance $\epsilon$,
the ID takes as input a matrix $A$, 
and returns the factorization $\wt{A} \Pi$ with
\begin{equation}\label{eq:idapprox}
\|A - \wt{A} \Pi \|_2 = O(\epsilon),
\end{equation}
where $\wt{A}$ is $n_t \times p$ and $\Pi$ is $p
\times n_f$. The parameter $p$ is 
the numerical rank of $A$ determined by the ID 
as part of the factorization. 
The columns of $\wt{A}$ are a $p$-column subset of the original
matrix $A$, chosen so that
the column space of $\wt{A}$ approximates that of $A$. The matrix
$\Pi$ contains the coefficients needed to approximately reconstruct the
columns of $A$ from those of $\wt{A}$. 
If we define the indexing function
$\pi$ so that the $m$th column of $\wt{A}$ is the $\pi(m)$th column
of $A$, then 
the approximation \eqref{eq:idapprox}
implies that
\[G(x_j^t,x_l^f) \approx \sum_{m=1}^p G(x_j^t,x_{\pi(m)}^f) \Pi_{ml}\]
for $l=1,\ldots,n_f$. Since the columns of $A$ represent the functions
$\{G(x,x_l^f)\}$ on a fine training grid, the expression above holds
not just for $x \in \{x_j^t\}$, but more generally for $x \in \Theta$.
That is,
\[G(x,x_l^f) \approx \sum_{m=1}^p G(x,x_{\pi(m)}^f) \Pi_{ml}.\]
Summing both sides of this expression
against $(\vec{\sigma})_l w_l^f$ and
rearranging yields
\[\sum_{l=1}^{n_f} G(x,x_l^f) (\vec{\sigma})_l w_l^f \approx \sum_{l=1}^{n_f}
  \sum_{m=1}^p G(x,x_{\pi(m)}^f) \Pi_{ml} (\vec{\sigma})_l w_l^f =
\sum_{m=1}^p G(x,x_{\pi(m)}^f) (\Pi W \vec{\sigma})_m\]
where $W$ is a diagonal $n_f \times n_f$ matrix with $W_{ll} = w_l^f$.
Since $\vec{\sigma} = B \hat{f}^K$, we let $T := \Pi W B$ to obtain
the representation \eqref{eq:outgoingcompressed} with
\begin{equation} \label{eq:Tmap}
  \vec{\rho} = T \hat{f}^K.
\end{equation}
$T$ is a generic $p \times K$
matrix which may be formed and stored 
once $\Pi$, $W$, and $B$ are available. We emphasize that each of these
matrices is identical for all patches of a given radius $\varepsilon$
and may therefore be precomputed.
$\Pi$ is obtained from a single interpolative decomposition, $W$ is
a simply a matrix of quadrature weights, and $B$ is computed by solving a
sequence of one-patch integral equations as explained in Section
\ref{sec:onepatchprelim}.

Using this compression scheme alone, it is straightforward to reduce the
cost of computing the sums \eqref{eq:mssums} from $\OO(K n_f N^2)$ to
$\OO(K p N^2)$.  The tools introduced in the remainder of this section
will allow us to reduce the cost further to $\OO(K p N \log N)$.

\subsection{Quadtree on the sphere}\label{sec:quadtreesphere}

\begin{figure}[ht]
  \centering
    \includegraphics[width=.35\linewidth]{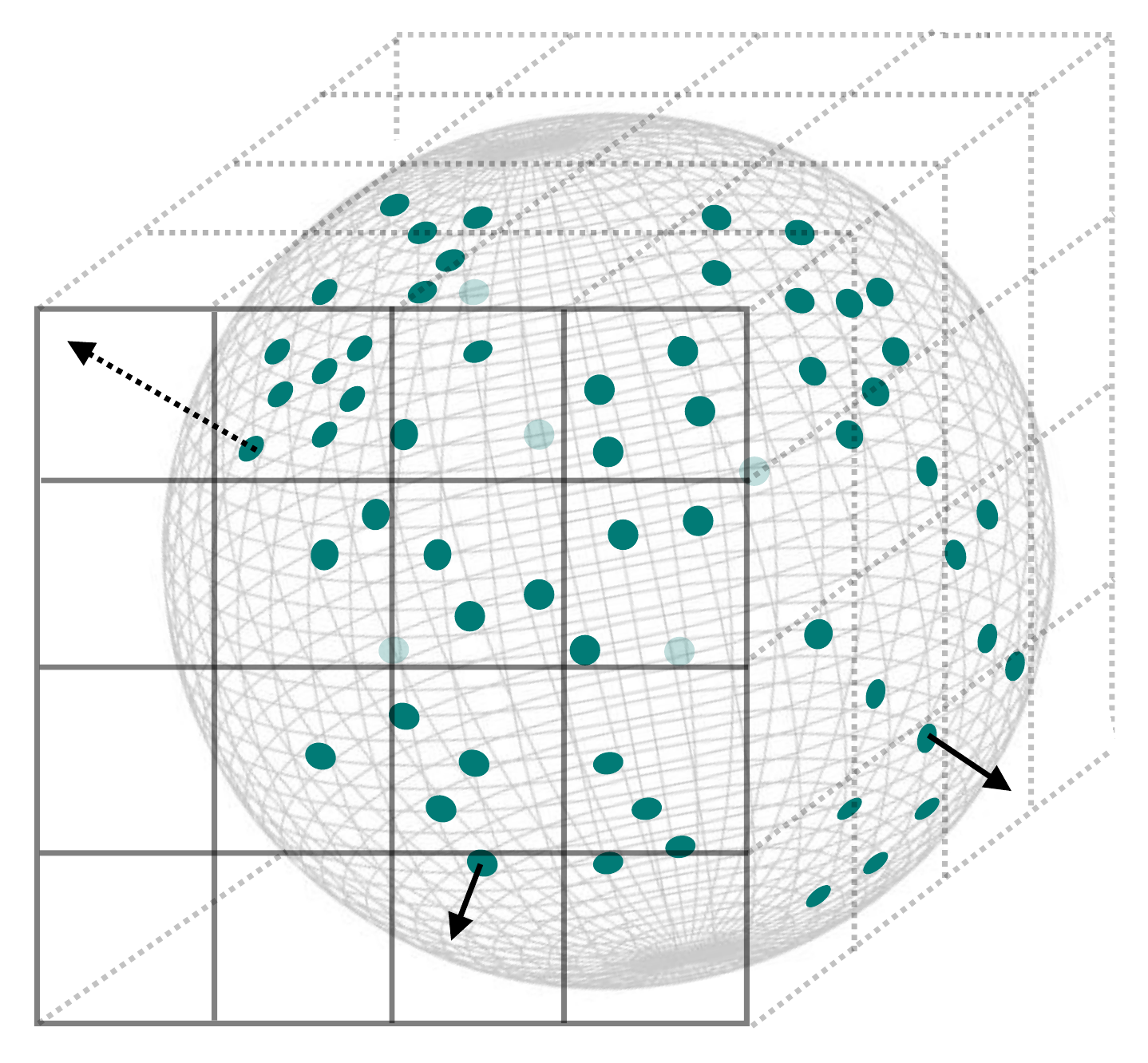} \hspace{.2in}
    \includegraphics[width=.32\linewidth]{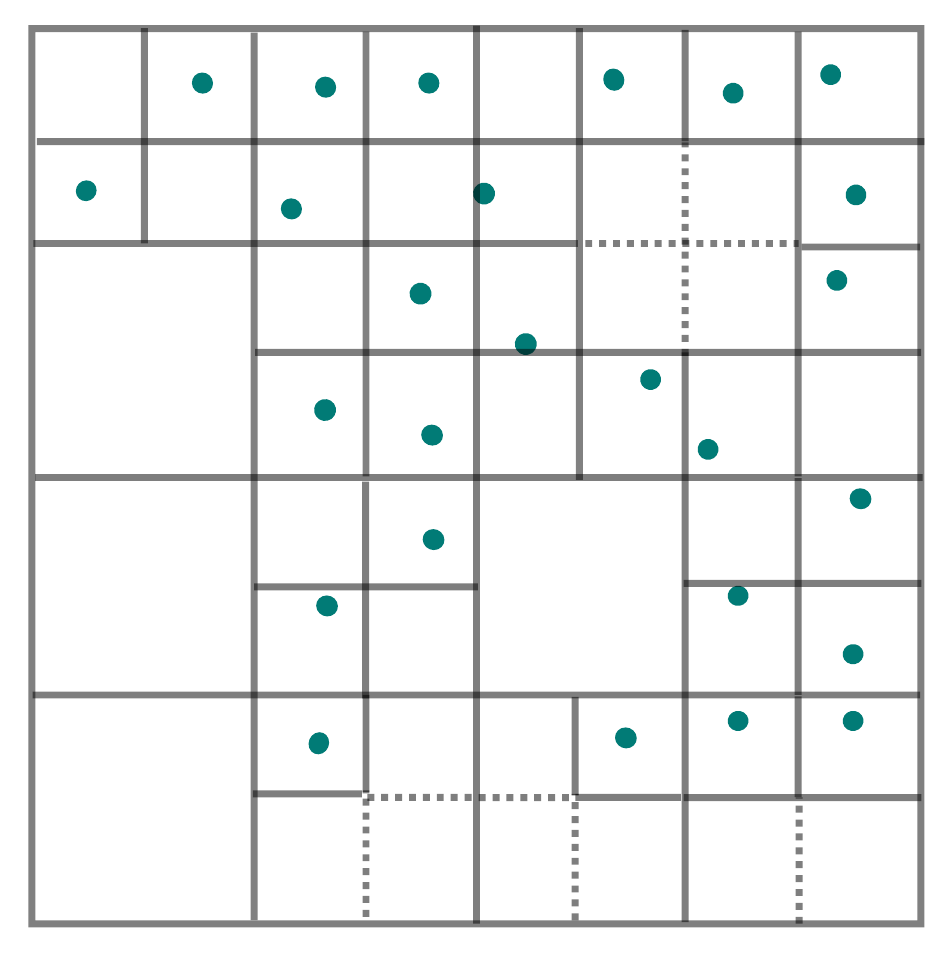}
    \caption{The sphere is inscribed in a cube and each 
    patch center is projected to a face of the cube by a ray
    emanating from the sphere center (left). An adaptive quad tree is
    then built on each face until, at the finest level, 
    there is one patch in every non-empty leaf node in the quad tree
    (right).}
    \label{fig:embeddedsphere}
\end{figure}

We now describe a data structure which will enable us to organize groups
of patches in a hierarchical fashion. We first inscribe the sphere in a
cube (see Fig. \ref{fig:embeddedsphere}).  We then project each patch
center onto the surface of the cube via the ray from the origin through
the patch center (indicated by the arrows in the figure).  This defines
a set of points on the surface of the cube. We then build a quadtree on
each face of the cube, subdividing boxes until there is only one point
per box, and pruning empty boxes in the process.  The union of these six
quadtrees is an FMM-like full tree data structure, which provides a
subdivision of the sphere itself into a hierarchy of levels.  The
patches assigned to a particular box in the full tree will be said to
form a {\em patch group}.  Each patch is a member of one patch group at
each level of the full tree. At the leaf level, each group consists of a
single patch.

We define parent, child, and neighbor boxes in the full tree 
in the same way as in an ordinary quadtree. The only modification to the
definition of a neighbor box is that it wraps across cube edges and corners.
Thus, a box adjacent to an edge has eight neighbors (like an interior box)
unless it is a corner box, in which case it has seven neighbors.
Well-separatedness and the interaction list for boxes or their 
corresponding patch groups are define as in the usual FMM. Two
boxes at a given level are well-separated if they are not
neighbors, and the interaction list for a particular box is comprised of
the well-separated children of its parent's neighbors.
We will sometimes refer to a patch $\Gamma_i$ as being in the
interaction list of some patch group $\gamma$, by which
we mean that $\Gamma_i$ is
contained in a group which is in the interaction list of $\gamma$.

\subsection{The representation of incoming fields on patch
groups}\label{sec:localrep}

\begin{figure}[ht]
  \centering
    \includegraphics[width=.4\linewidth]{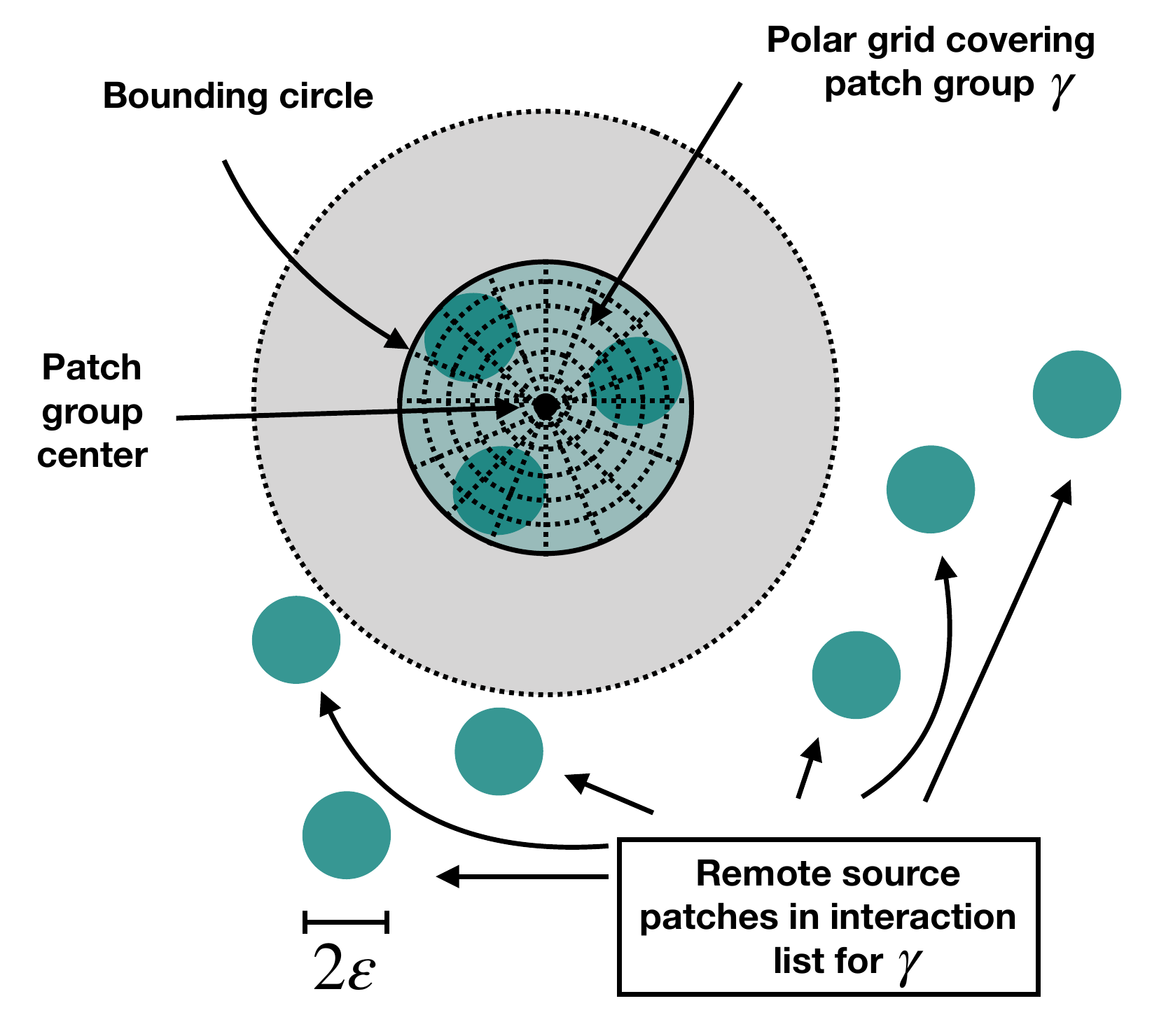}
    \caption{
    For a group of $m$ patches, the field due to well-separated 
    source patches may be captured with high order accuracy on a polar
    grid which covers all $m$ patches. 
    }
    \label{fig:localgrids}
\end{figure}

Since the incoming field due to remote source patches in the interaction
list of a patch group $\gamma$ is smooth, it can be efficiently
represented on a spectral polar grid (see Fig.
\ref{fig:localgrids}). This requires the construction of a {\em bounding
circle} on the surface of the sphere, enclosing all of the patches in
$\gamma$, which circumscribes the grid. Incoming field values can
then be obtained at arbitrary points inside the bounding circle by interpolation. We refer to the
grid samples of the incoming field as an {\em incoming representation}.

The bounding circle is straightforward to construct using a 
``smallest circle algorithm" for a collection of 
points in the plane, suitably adapted to the sphere (see 
\cite{skyum91,welzl91,xu03} and the references therein for discussion of
the smallest circle problem).

Given a bounding circle for a patch group, we can build a local polar
coordinate system $(r,\theta)$, for which $r = 0$ corresponds to the
center of the patch group, and $r = R$ corresponds to the bounding
circle. We must select an {\em incoming grid} in these coordinates 
which can represent a
smooth incoming field in a high order manner with as few grid points
as possible.  For this, we will use
a parity-restricted Chebyshev-Fourier
basis, formed by
taking products of scaled Chebyshev polynomials in the radial variable 
$r \in [-R,R]$ with
trigonometric functions in the angular variable 
$\theta \in [0, 2 \pi)$. The
coefficients of an expansion in these basis functions
corresponding to Chebyshev and Fourier modes of different parity can be
shown to be zero, hence the name of the basis.
This is an efficient and spectrally accurate basis with a simple
associated grid \cite{boyd11}. 
Namely, the coefficients of the expansion may be computed from function samples
 on a polar grid comprised of the scaled Chebyshev nodes 
in $r \in [0,R]$
and equispaced nodes in $\theta \in [0,2 \pi)$. 
The desired field may then be evaluated at any point
inside a patch group's bounding circle
by evaluating the resulting Chebyshev-Fourier expansion.
It is straightforward to verify that the number of grid points 
and coefficients required to obtain an accuracy $\epsilon$ is $\OO(\log^2(1/\epsilon))$.

\section{Solution of the multiple scattering system}\label{sec:algorithm}

We now describe our method to solve the discretized many-patch system
\eqref{eq:opinteqnsdisc}, including the fast algorithm for accelerating
the computation of the multiple scattering interactions
\eqref{eq:mssums} within a GMRES iteration. 

\vspace{.2in}

\noindent {\bf Step 1: Precomputation (for each choice of $\varepsilon$)}

Given the patch radius $\varepsilon$, select
the Zernike truncation parameter $K$ and form
the matrix $Q$. 

(a) Solve the system $S B = Q$ described in Section \ref{sec:onepatch}.

(b)
Construct the matrix $T$ defined in Section
    \ref{sec:ID} by building and composing the matrices $\Pi$,
    $W$, and $B$. $\Pi$ need not be stored after $T$ is formed.
  
(c)
Construct the vector 
$I = (w_1^f,\ldots,w_{n_f}^f) B$, 
used to obtain the quantities $J$ and $\mu$ in \eqref{eq:getI}.
At this point we no longer need to store $B$,
    only the $p \times K$ matrix $T$ and the $1 \times K$ vector $I$. 
The storage associated with the outputs of the precomputation phase is
therefore negligible.

\vspace{.2in}

\noindent {\bf Step 2: Construction of hierarchical data structure}

Let $N$ denote the number of patches on the surface of the sphere,
assumed to satisfy the
the minimum patch separation
condition introduced in Section \ref{sec:manypatchsetuptools}.

(a) Form the quadtree on the sphere described in Section
    \ref{sec:quadtreesphere}. The data structure should associate
    each patch with its group at every level, and identify
    the interaction list of every patch group.

(b) For each patch group, construct the incoming grid described in
Section \ref{sec:localrep}. For each patch, construct the Zernike
sampling grid described in Section \ref{sec:zernike}.

\vspace{.2in}

\noindent {\bf Step 3: Iteration}

We use GMRES to solve the system \eqref{eq:opinteqnsdisc}. At
each iteration, we must apply the system matrix; that is, we must compute
\begin{equation}\label{eq:sysmatapply}
  \hat{f}_i^K + P \sum_{j \neq i} S_{ij} B \hat{f}_j^K 
\end{equation}
for $i=1,\ldots,N$, where here $(\hat{f}_1^K,\ldots,\hat{f}_N^K)^T \in
\RR^{KN}$ is the input vector at a given iteration. 
The following algorithm computes this
expression in $\OO(N \log N)$ operations. 

\begin{enumerate}
  \item Compute and store the outgoing
    coefficients $\vec{\rho}_i = T \hat{f}_i^K$ for each patch, $i=1,\ldots,N$.

    {\em \underline{Cost}: Approximately $p K N$.}

  \item Loop through every patch group in every level. For each patch group $\gamma$, loop through all patches in its 
    interaction list. For each such patch $\Gamma_i$, evaluate the
    field induced by the density on $\Gamma_i$ on the
    incoming grid of $\gamma$, using the outgoing representation
    \eqref{eq:outgoingcompressed}. Add together all such field values to obtain the total incoming
    field on the incoming grid.

    {\em \underline{Cost}: If $q$ is an upper bound on the number of points in each
    incoming grid, the cost of evaluating a single outgoing representation on
    an incoming grid is at most $q p$. At each level, the outgoing
    representation corresponding to each patch must be
    evaluated on at most $27$ incoming grids, since the
    interaction list of each patch's group at that level contains at
    most $27$ other groups. There are approximately $\log_4 N$ levels.
    Therefore, the
    cost of this step is approximately $27 q p N \log_4 N$.}

  \item At the leaf level of the tree, each patch group $\gamma$ contains a
    single patch, say $\Gamma_i$. Though we have already evaluated the outgoing
    representation for $\Gamma_i$ on the incoming grids of all
    (single-patch) groups in the interaction list of $\gamma$, we now do
    so also for the neighbors of $\gamma$, which are
    also single-patch groups but are not contained in the interaction
    list of $\gamma$. We add these contributions to the field values
    already stored on the incoming grids of these neighbor
    patches.

    {\em \underline{Cost}: Since each leaf-level single-patch group has at most $8$ neighbors, 
    the cost of this step is approximately $8 q p N$.}

    {\em Note: For each patch $\Gamma_i$, the incoming field due to every other
    patch has now been stored in the incoming grid of exactly one patch-group
    of which $\Gamma_i$ is a member. Indeed, every other patch is
    either a neighbor of $\Gamma_i$ at the leaf level, or it is contained
    in exactly one of the interaction lists of the patch groups
    containing $\Gamma_i$.}

  \item Loop through each patch group. For every patch $\Gamma_i$ in a
    group $\gamma$, evaluate the interpolant of the incoming field stored
    on the incoming grid of $\gamma$ at the Zernike sampling nodes
    on $\Gamma_i$.

    {\em \underline{Cost}: There are $\OO(K)$ Zernike sampling nodes,
    so the cost of each
    interpolation is approximately $q^2$ to form the interpolant and $K
    q$ to evaluate it. Each patch is a member of a single group at each
    level, so we must carry out approximately $N \log_4 N$ such
    interpolations. The total cost is therefore approximately $(q^2
    + Kq) N \log_4 N$. (For large $q$, this step could be
    accelerated with fast transform methods but $q$ is generally too small
    for this to provide any significant benefit.)}

    At this point, we have computed the field due to all other patches
    on the Zernike sampling grid on each patch. That is, we have
    computed the sums $\sum_{j \neq i} S_{ij} B \hat{\sigma}_j$ for $i =
    1,\ldots,N$.

  \item Apply the matrix $P$ to the values stored on the Zernike
    sampling grid on each patch and add $\hat{f}_i^K$ to the result to
    obtain \eqref{eq:sysmatapply}.

    {\em \underline{Cost}: Approximately $K^2 N$.}
\end{enumerate}

The total cost of each iteration is therefore $\OO(N \log N)$, with
asymptotic constants which involve the parameters $K$, $q$, and $p$
associated with the resolution of smooth functions on spectral grids.
The singular character of the problem is dealt with entirely during the
precomputation phase.

\subsection{Optimizations and parallelization}

While the algorithm described above has the desired computational
complexity, there are several practical considerations that are worth
discussing to optimize its performance.

{\em Selection of incoming grid parameters:}
Rather than making a uniform choice of the radial and azimuthal
truncation parameters for the incoming grid, we can compute these
adaptively as follows. For each patch group $\gamma$, we determine the
distance from its bounding circle to the nearest patch in its
interaction list. We then adaptively construct an incoming grid which
accurately interpolates a collection of point sources $G(x,x')$ at
points $x'$ this distance away. This adaptive interpolation is carried
out by increasing the incoming grid truncation parameters until the last
few Legendre-Fourier coefficients of the interpolant fall below some
specified tolerance.

{\em Additional compression of the outgoing representation:}
Instead of using the same outgoing coefficients $\vec{\rho}_i$ for each
level of the quadtree, we can associate with each patch a different
outgoing representation for each level. Recall that the far field
regions $\Theta_i$ were constructed
identically for each patch $\Gamma_i$ to be as large as possible,
consistent with the minimum patch separation. This way, one could build
a single generic matrix $T$ taking a density on a patch to its outgoing
representation. $T$ was built by compressing the outgoing field due to a
generic patch $\Gamma$ against a grid on a generic far field region $\Theta$.
Instead, we can build one such matrix for each level of the quadtree by
constructing a generic far field region for each level. Each such far field
region is an annulus or disk on the surface of the sphere. For each
level, it is taken to be just large enough so that for any $i =
1,\ldots,N$, in the coordinate system of $\Gamma_i$, it covers the
bounding circle of every group $\gamma$ containing $\Gamma_i$ in its
interaction list at that level.  Using the interpolative decomposition,
we can then recompress the outgoing representation for a generic patch
against training grids on each of the approximately $\log_4 N$ new far
field regions. We obtain one matrix $T$ per level, each of which has
fewer rows and therefore yields fewer outgoing coefficients than the
original.

{\em Parallelization:}
Each step of the algorithm to compute \eqref{eq:sysmatapply} may be 
straightforwardly
parallelized. Steps (1) and (5) are parallelized over all patches; steps
(2) and (4) are parallelized over all patch groups at all levels; step
(3) is parallelized over all patch groups at the leaf level. 

\section{The one-patch
integral equation}\label{sec:onepatch}

In this section, we describe in detail a solver for the
integral equation \eqref{eq:onepatinteq}, 
as well as the construction of the far-field quadrature nodes
$x_{i,1}^f,\ldots,x_{i,n_f}^f$ and weights
$w_1^f,\ldots,w_{n_f}^f$ discussed in Section \ref{sec:onepatchprelim}.

We assume that a patch $\Gamma$ has radius
$\varepsilon$ and make use of
cylindrical coordinates $(r,\theta,z)$. If we take the center of the
patch to be the north pole of the sphere, then $r = 0$ corresponds to
the $z$-axis, $r = 0$ and $z = \pm 1$ to the north and south poles, respectively, and $\theta =
0$ to the $x$-axis.
Following the approach of \cite{young12,helsing14}, we use the
rotational symmetry of $\Gamma$ to reduce the integral
equation over the patch to a sequence of one-dimensional integral equations, 
each corresponding to a
Fourier mode in the variable $\theta$. 
More precisely, 
we denote by $C$ the arc which generates $\Gamma$ via rotation
about the $z$-axis:
$C(t) \equiv (r(t),z(t)) =
(\sin(t),\cos(t))$ for $t \in [0,\varepsilon]$. 
In this parametrization, $t$ is simply the arclength along the sphere.

Let $x = (r,\theta,z)$ and $x' = (r',\theta',z')$. 
Since $G_E$ and $G_I$ are functions of $|x-x'|$ 
and
\[|x-x'| = \sqrt{r^2 + r'^2 + (z-z')^2 -
2rr'\cos(\theta-\theta')},\]
we can write the dependence of the Green's function in cylindrical coordinates as
$G(x-x') = G(r,r',z-z',\theta-\theta')$. 
In these coordinates, the one-patch integral equation
\eqref{eq:onepatinteq} takes the form 
\[\int_0^\varepsilon 
\int_0^{2 \pi} G(r(t),r'(t'),z(t)-z'(t'),\theta-\theta') 
\sigma(r'(t'),z'(t'),\theta') r'(t') \, dt' \, d \theta' = f(r(t),z(t),\theta).\]
Representing $\sigma$ as a Fourier series in $\theta$,
\[\sigma(r(t),z(t),\theta) = \sum_{n=-\infty}^\infty \sigma_n(t) e^{i n
\theta},\]
and taking the Fourier transform of both sides of this equation, upon
rearrangement, gives the following integral equation for the Fourier
modes:
\begin{equation}\label{eq:modalinteqn}
  2 \pi \int_0^\varepsilon G_n(t,t') 
\sigma_n(t') \sin(t') \, dt' = f_n(t).
\end{equation}
Here $G_n(t,t')$, $\sigma_n(t)$,
and $f_n(t)$ are the Fourier transforms of
$G(r(t),r'(t'),z(t)-z'(t'),\theta)$, 
$\sigma(r(t),z(t),\theta)$ and 
$f(r(t),z(t),\theta)$ with respect to $\theta$.
Thus, after solving the one-dimensional modal equations 
\eqref{eq:modalinteqn},
we can recover
$\sigma(r(t),z(t),\theta)$ from its Fourier series. Note that
the Fourier series is spectrally convergent because
$\sigma(r(t),z(t),\theta)$ is smooth as a function of $\theta$, even though 
it is singular as a function of $t$ at the edge $t = \varepsilon$.

\subsection{Evaluation of the modal kernels}

Let
\begin{align*}
  G_n^{(1)}(t,t') &= \frac{1}{\pi} \int_0^\pi \frac{2}{|x-x'|} 
  \cos(n \tilde{\theta}) \, d \tilde{\theta} \\
  G_n^{(2)}(t,t') &= \frac{1}{\pi} \int_0^\pi
\log\left(\frac{2}{|x-x'|} \right) \cos(n \tilde{\theta}) \, d \tilde{\theta} \\
  G_n^{(3)}(t,t') &= \frac{1}{\pi} \int_0^\pi \log \left( 1 +
\frac12 |x-x'| \right) \cos(n \tilde{\theta}) \, d \tilde{\theta}.
\end{align*}
Then, using the formulae \eqref{eq:Gextonsurface} and \eqref{eq:Gintonsurface},
it is straightforward to show that
$G_n = G_n^{(1)} + G_n^{(2)} - G_n^{(3)}$ 
for $G_E(x,x')$ and 
$G_n = G_n^{(1)} - G_n^{(2)} - G_n^{(3)}$ 
for $G_I(x,x')$.
We can write $|x-x'|$ in terms of $t$, $t'$ and
$\tilde{\theta} = \theta-\theta'$
as
\[|x-x'| = \sqrt{2 \left( 1-\cos(t)\cos(t') -
\sin(t)\sin(t')\cos(\tilde{\theta}) \right)}.\]
The integrands are not smooth at $t = t'$, $\tilde{\theta} = 0$, so we
must use specialized methods to evaluate each kernel.

$G_n^{(1)}(t,t')$ is simply the cosine transform of the Coulomb kernel
and arises in boundary integral equations for
electrostatics on axisymmetric surfaces. In 
\cite{helsing14}, an efficient evaluation algorithm is described
which involves writing the modal kernel in terms of
Legendre functions of half-integer order and using their
associated three-term recurrence.  We refer the reader to this 
paper for further details.

The kernel $G_n^{(2)}(t,t')$ is weakly singular and may be evaluated by
adaptive Gaussian quadrature. However, the following formula, discovered
by a combination of analytical manipulation and symbolic calculation
with Mathematica, has been numerically verified for a wide
range of values and is significantly faster:
\begin{equation*}
\frac{1}{\pi}\int_0^\pi \log\left( \frac{2}{|x-x'|} \right) \cos(n
\tilde{\theta}) \, d \tilde{\theta} = 
\begin{cases}
  -\log\left(\cos(t_1/2)\sin(t_2/2)\right) & n = 0 \\
  \frac{1}{2 n} \left(\tan(t_1/2) \cot(t_2/2)\right)^n & n > 0 \\
  t_1 = \min(t,t'), t_2 = \max(t,t').
\end{cases}
\end{equation*}

The integrand in the expression for $G_n^{(3)}(t,t')$ is even more weakly
singular, so $G_n^{(3)}(t,t^\prime)$ may be evaluated relatively quickly
by adaptive Gaussian quadrature.

\subsection{Discretization of the modal integral equations}

Since \eqref{eq:modalinteqn} is a singular integral equation, care must
be taken to discretize it accurately. The dominant singularity of the kernel $G_n(t,t')$ at $t = t'$ 
is the logarithmic singularity of $G_n^{(1)}(t,t')$.
An analogous classical problem is therefore the first-kind integral equation
arising from the solution of the Dirichlet problem on an open arc in two
dimensions by a single layer potential.
Stable and accurate numerical schemes for this problem can be found, for example,
in \cite{yan88,atkinson91,jiang04}. As described in 
\cite{jiang04}, when the domain is the
interval $[-1,1]$, the solution of
\begin{equation}
\label{straighteq} 
\int_{-1}^1 \log|t-s| \sigma(s) \, ds = f(t)
\end{equation}
can be computed with spectral accuracy in the form 
$\sigma(t) = g(t)/\sqrt{(1+t)(1-t)}$,
where $g$ is a smooth function whose Chebyshev coefficients depend in a
simple manner on those of $f$.
For an open
arc, the corresponding integral equation can be preconditioned
using the solution of \eqref{straighteq}. This procedure results in
a Fredholm equation of the second kind 
for which the density may be represented as a Chebyshev expansion
and computed stably with high order accuracy.

In the present context, the inclusion of the additional 
weakly singular kernels $G_n^{(2)}$ and $G_n^{(3)}$ 
cause the singularity of $\sigma_n(t)$ to be more complex, but our
numerical evidence suggests that there is still a dominant square root
singularity at $t = \varepsilon$. To be more precise, if we represent
$\sigma_n$ by 
\begin{equation}
\sigma_n(t) = g_n(t)/\sqrt{\varepsilon - t}
\label{sigmangndef}
\end{equation}
near $t = \varepsilon$, we can investigate the effectiveness of
representing $g_n$ in a basis of orthogonal polynomials.
While the exact behavior of $g_n(t)$
is not understood analytically, the numerical results presented in Section
\ref{sec:sigsingularity} suggest that it is only
mildly non-smooth. We note that there is no singularity at the endpoint
$t = 0$, since this point corresponds to the patch center, at which
there is no physical singularity.

To resolve the endpoint singularity of $\sigma_n$, we discretize
it on a set of panels $[a_0,a_1],[a_1,a_2],\ldots,[a_{m-1},a_m]$ on
$[0,\varepsilon]$ which are dyadically
refined towards $t = \varepsilon$:
\[ a_0 = 0,\ a_1 = \frac{\varepsilon}{2},\ a_2 = \frac{3 \varepsilon}{4},
\dots,\ a_{m-1} = \frac{(2^{m-1}-1) \varepsilon}{2^{m-1}},\ a_{m} =
\varepsilon.\]
On each panel, except the last, $\sigma_n$ is represented as a Legendre series of fixed
order $k$. Since $\sigma_n$ is smooth on each such
panel and separated
from its singularity by a distance equal to the panel length, it can be
shown that
this representation has an
error of size $\OO(e^{-k} \log_2 (1/\varepsilon))$.
This argument is widely used in handling endpoint and corner singularities in the context
of boundary integral equations 
\cite{bremer1,bremer2,helsing1,helsing3,serkh2016jcp,trefethen}.

On the last panel, we analytically incorporate a square root singularity
into our representation of $\sigma_n$ as above, and expand $g_n(t) =
\sigma_n(t) \sqrt{\varepsilon - t}$ as a
series of Jacobi polynomials with $\alpha = -\frac12$ and $\beta = 0$.
If the singularity of $\sigma_n$ at $t = \varepsilon$ were exactly of
square root type, this would yield a spectrally accurate representation
of $\sigma_n$. Instead, as we will show in Section
\ref{sec:sigsingularity}, we obtain a
representation which is finite order but resolves the solution quite
well even for modest truncation parameters.

Thus we have rewritten \eqref{eq:modalinteqn} as
\[
  f_n(t) 
= 2 \pi \sum_{j=1}^{m-1} \int_{a_{j-1}}^{a_j} G_n(t,t')
\sigma_n(t') \sin(t') \, dt' \\ 
+ 2 \pi \int_{a_{m-1}}^{\varepsilon}
\frac{G_n(t,t')}{\sqrt{\varepsilon - t'}} \left( \sigma_n(t')
\sqrt{\varepsilon - t'} \right) \sin(t') \, dt' \, 
\]
and discretized $\sigma_n$ by Legendre polynomials for the first $m-1$
panels and by Jacobi polynomials for the last. Sampling the resulting
equations at the corresponding quadrature nodes - Gauss-Legendre for the
first $m-1$ panels and Gauss-Jacobi for the last - yields a collocation method for
$\sigma_n$, in which $\sigma_n$ is determined by its piecewise polynomial
basis coefficients. For each collocation node $t_i$, we compute the
system matrix entries by
adaptively integrating $G_n(t_i,t')$ in $t'$ against
the piecewise polynomial basis functions. We compute
the values $f_n(t_i)$ by discretizing the Fourier
transform of $f(r(t_i),z(t_i),\theta)$ in $\theta$ by the trapezoidal
rule, which is spectrally accurate for smooth, periodic functions. We
solve the resulting set of linear systems - one for each Fourier mode - by $LU$ factorization and
back substitution. The factorizations may be reused, since we must solve
a one-patch integral equation for many different right hand sides.

We can now define the fine grid points and the smooth quadrature weights
introduced in Section \ref{sec:onepatchprelim}. The points
$x_{i,1}^f,\ldots,x_{i,n_f}^f$ are the tensor products of the
collocation nodes in the radial direction with equispaced points - the
trapezoidal rule quadrature nodes - in the azimuthal direction.
$w_1^f,\ldots,w_{n_f}^f$ are the corresponding quadrature weights -
products of the panel-wise Gauss weights with the trapezoidal rule weight.

\subsection{Numerical investigation of the singularity of
$\sigma_n$}\label{sec:sigsingularity}

In this section, we contrast two strategies for representing $\sigma_n$
in \eqref{sigmangndef}.
In the first, we use $m = 1$ panels, and represent
$g_n$ in a basis of Jacobi polynomials, which takes into account the square
root singularity in $\sigma_n$. This approach would yield spectral
accuracy with respect to $g_n$ if $\sigma_n$ only contained a square
root singularity. The second strategy is the one described above; we
use $m > 1$ panels with a Jacobi polynomial basis of fixed degree only in the
last panel. These experiments give us some insight into the nature of the
true singularity in $\sigma_n$, and justify our discretization choice.

In both cases, we solve the interior one-patch integral equation by the method
described above for a basis of Zernike polynomials with truncation
parameter $M = 15$. The results do not change significantly if we solve
the exterior equation instead. We do this for several different choices of
$\varepsilon$. The Fourier series truncation is fixed sufficiently large
to resolve the highest azimuthal Zernike mode. For each solution, 
we measure the residual error in $L^2$, normalized by the patch size:
\begin{equation}\label{eq:patcherrl2}
  \norm{\calS \sigma - f}_{L^2(\Gamma)} / |\Gamma|.
\end{equation}
Here $|\Gamma|$ is the surface area of the patch, and $f$ is
a Zernike polynomial. 
This measures the extent to which the
computed solution of the one-patch BVP satisfies the Dirichlet boundary
condition. This solution automatically satisfies the Neumann boundary
condition and the PDE, because of its representation as a single layer
potential with the Neumann Green's function, so a small $L^2$ residual error
corresponds to a solution which nearly satisfies the boundary value
problem. This error is computed by
quadrature on a Legendre-Fourier grid which does not overlap with the
grid on which the integral equation is solved, so it is {\em not} the same as
the residual of the solution to the discrete linear system.

\begin{figure}[ht]
  \centering
    \includegraphics[width=.45\linewidth]{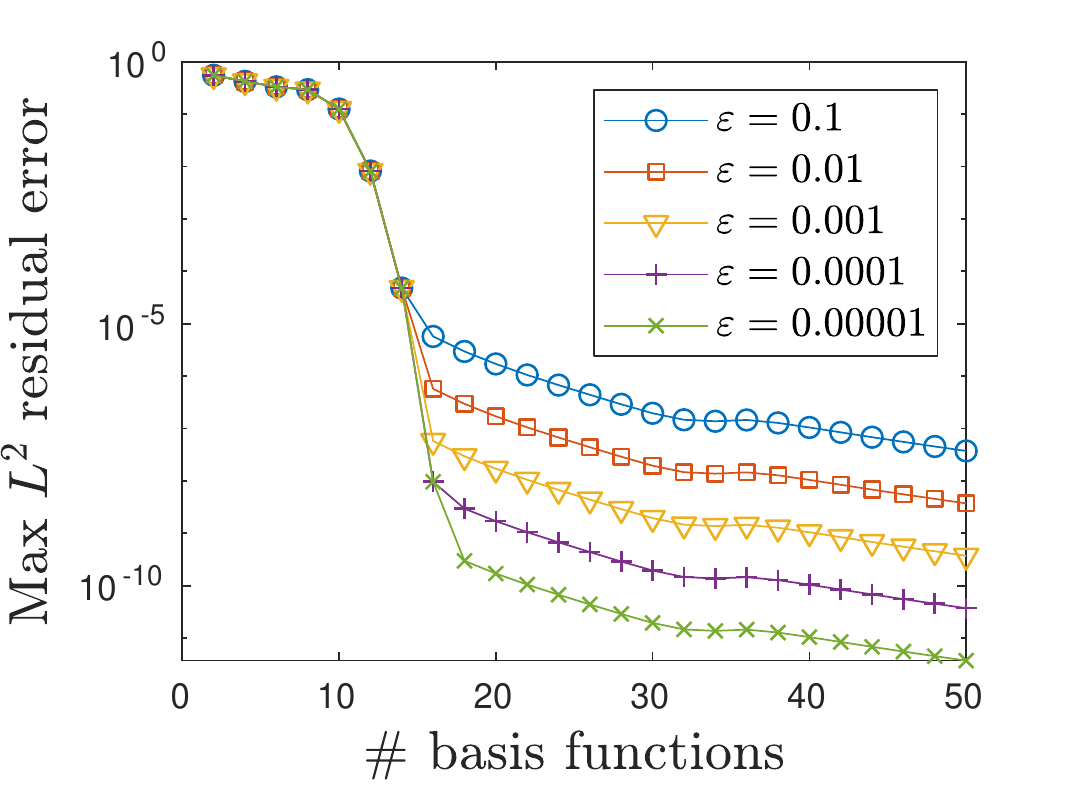} \hspace{.2in}
    \includegraphics[width=.45\linewidth]{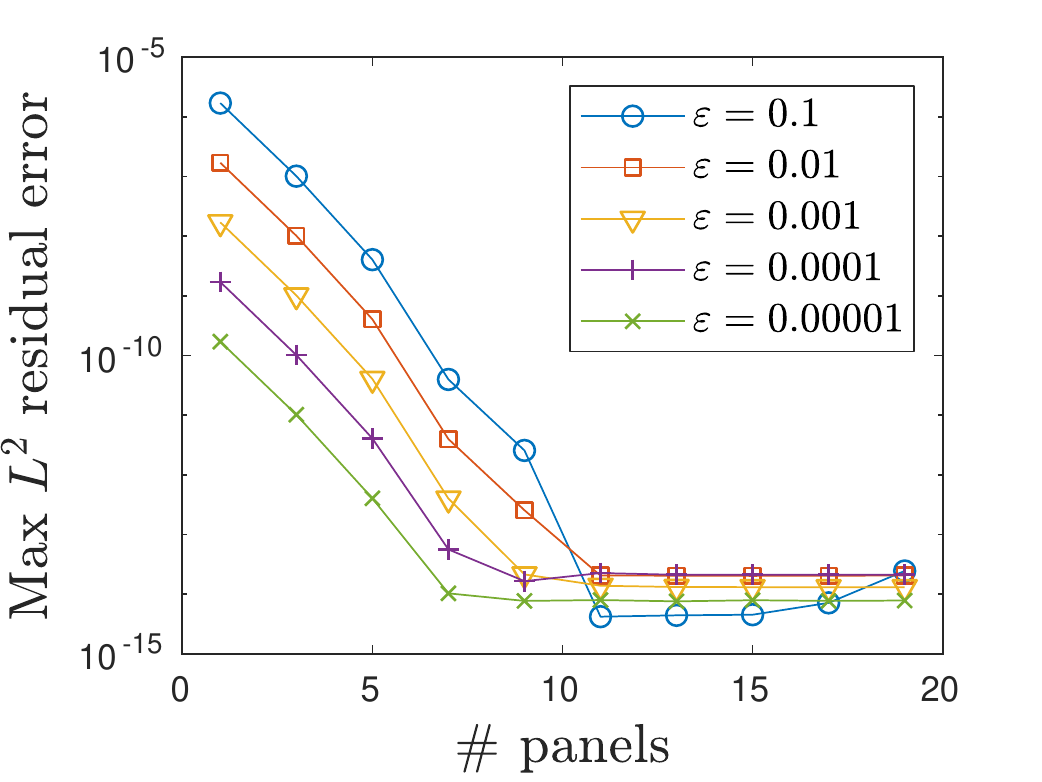}
    \caption{Left panel: $g_n$ is represented by a basis of 
    Jacobi polynomials on a single panel. We plot the
    maximum residual error \eqref{eq:patcherrl2} vs.
  the number of Jacobi basis functions. Right panel: $g_n$ is represented in
a Legendre basis on every panel except the last, where a Jacobi basis is used.
We plot the maximum residual error vs. the number of panels.}
    \label{fig:onepaterr}
\end{figure}

Using the first strategy ($m=1$), we 
measure the error
\eqref{eq:patcherrl2} for each Zernike polynomial, as the number of
Jacobi basis functions is increased. The error is defined to be the
maximum taken over all Zernike polynomials.
The results are presented in the
left panel of Fig. \ref{fig:onepaterr}. We observe an initial regime of rapid convergence, followed by
much slower convergence. Indeed, $15$ basis functions are required
to resolve the highest Zernike modes we have used as data. Afterward, the slow
regime of convergence suggests that $\sigma_n$ has a dominant square root singularity
and a subdominant term which is nonsmooth, but much smaller. We also notice
that performance improves as $\varepsilon$ is decreased, which is not
surprising since as $\varepsilon \to 0$, we approach the flat case in
which $\sigma_n$ has a pure square root singularity.

The second strategy is explored in the right panel of Fig
\ref{fig:onepaterr}. Here, we fix $20$
basis functions per panel - sufficient to begin with a good error
constant, according to the first experiment. We then increase the number 
$m$ of panels. Although we can already obtain
quite good accuracy using the first strategy, the second allows us to
reach near-machine precision. The improvement is particularly dramatic
for larger choices of $\varepsilon$.

\section{Numerical experiments}\label{sec:numresults}

An important parameter in studying narrow escape and narrow capture problems
is the {\em patch area fraction} $f_{N,\varepsilon}$. 
Since the surface area
of a single patch of radius $\varepsilon$ is given by
\[A_\varepsilon = 4 \pi \sin^2(\varepsilon/2), \]
we have
\begin{equation}
  f_{N,\varepsilon} = N \sin^2(\varepsilon/2).
\end{equation}
Assuming $\varepsilon$ is sufficiently small, we may write
\begin{equation}\label{eq:patareafrac}
  f_{N,\varepsilon} \approx \varepsilon^2 N / 4.
\end{equation}
Given $N$, we will use \eqref{eq:patareafrac} to compute
the patch radius $\varepsilon$ for a given patch area fraction.

\subsection{Convergence with respect to the Zernike basis} \label{sec:Mconv}

We first investigate the convergence
of the solution with respect to the Zernike truncation parameter $M$,
which determines the largest radial and azimuthal Zernike modes used to
represent the smooth incoming field on each patch. 
We fix the patch area fraction at $f_{N,\varepsilon} = 0.05$ and
carry out 
experiments with $N = 10$, $100$, and $1000$ patches. $\varepsilon$
is computed from \eqref{eq:patareafrac}. 
The patch locations are drawn from a uniform random
distribution on the sphere, with a minimal patch separation of $2
\varepsilon$ enforced. 
In each case, we solve the one-patch problems with the
truncation parameter $M$ set to $1,3,5,\ldots,15$.
The one-patch solutions are obtained, guided
by the results in Fig. \ref{fig:onepaterr}, 
using $13$ panels with $20$ basis functions
per panel, and the number of Fourier modes set equal to the number of
azimuthal modes in the Zernike basis. The ID and GMRES tolerances are
set to $10^{-15}$, and the incoming grid 
tolerance is set to $10^{-12}$. 

We measure error in two ways. The first,
as in \eqref{eq:patcherrl2}, is to examine the relative $L^2$ residual
of the multiple scattering system \eqref{eq:opinteqns}
(the discrepancy of the computed boundary values with the Dirichlet data) 
on a random patch $\Gamma_i$:
\begin{equation}\label{eq:l2res}
  \frac{1}{|\Gamma_i|} \norm{\left(\calS \sigma_i + \sum_{j \neq i}^N \calS_{ij} \sigma_j
  \right) - 1}_{L^2(\Gamma_i)}.
\end{equation}
The second is to examine the difference between the computed average 
mean first passage time (MFPT) $\mu$ and
a reference value, denoted by $\mu_{\text{ref}}$. 
We obtain $\mu_{\text{ref}}$ by carrying out a more refined simulation,
with $M = 17$ on each patch, while also increasing
the number of panels and basis
functions used to solve the one-patch problem to $19$ and $30$,
respectively, and doubling the numbers of both radial and azimuthal modes 
used in the incoming grids of all patch groups. 
This is a self-consistent convergence
test for $\mu$.

\begin{figure}[ht]
  \centering
    \includegraphics[width=.32\linewidth]{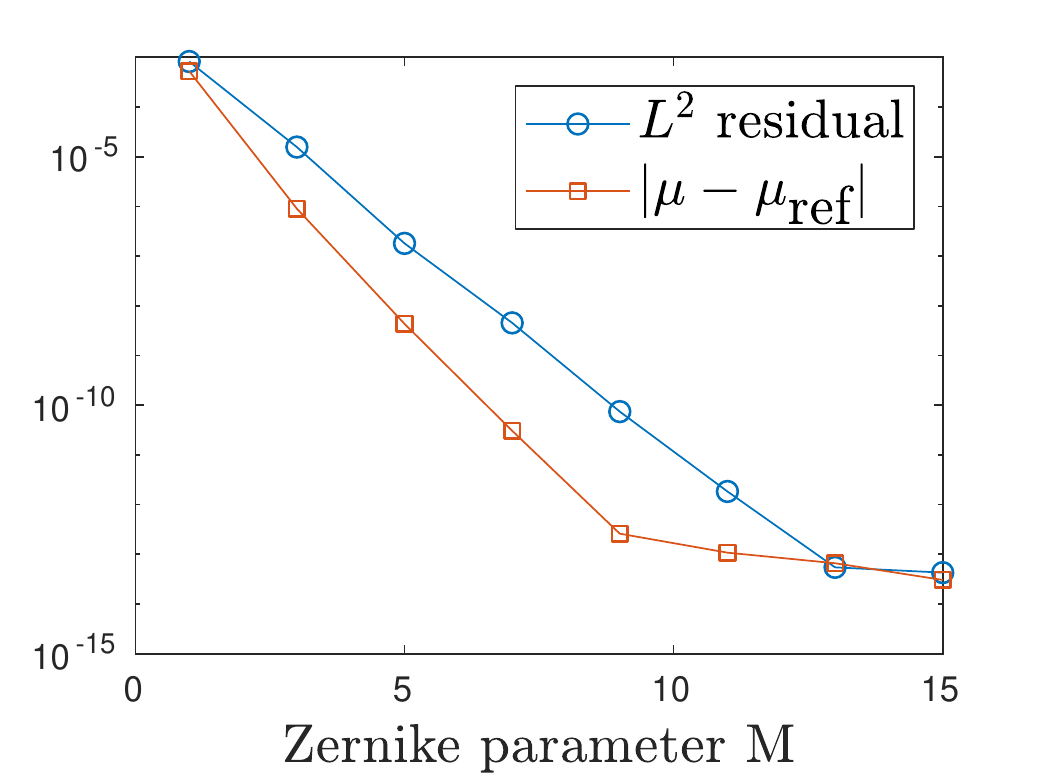} \hspace{.01in}
    \includegraphics[width=.32\linewidth]{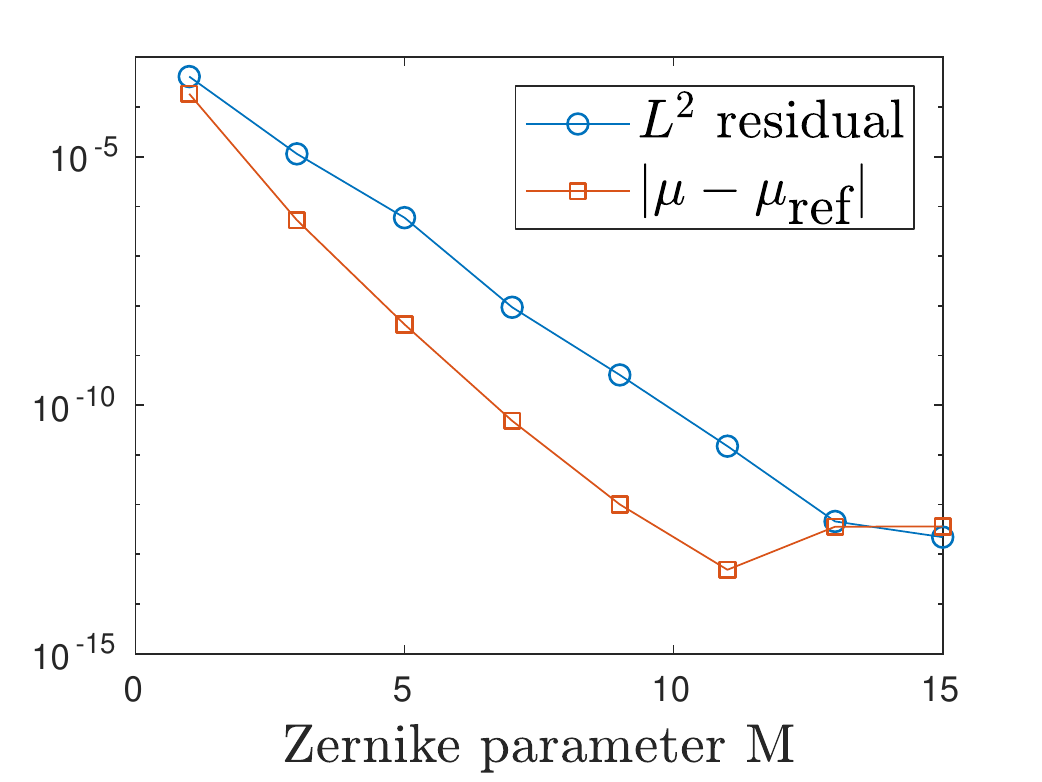} \hspace{.01in}
    \includegraphics[width=.32\linewidth]{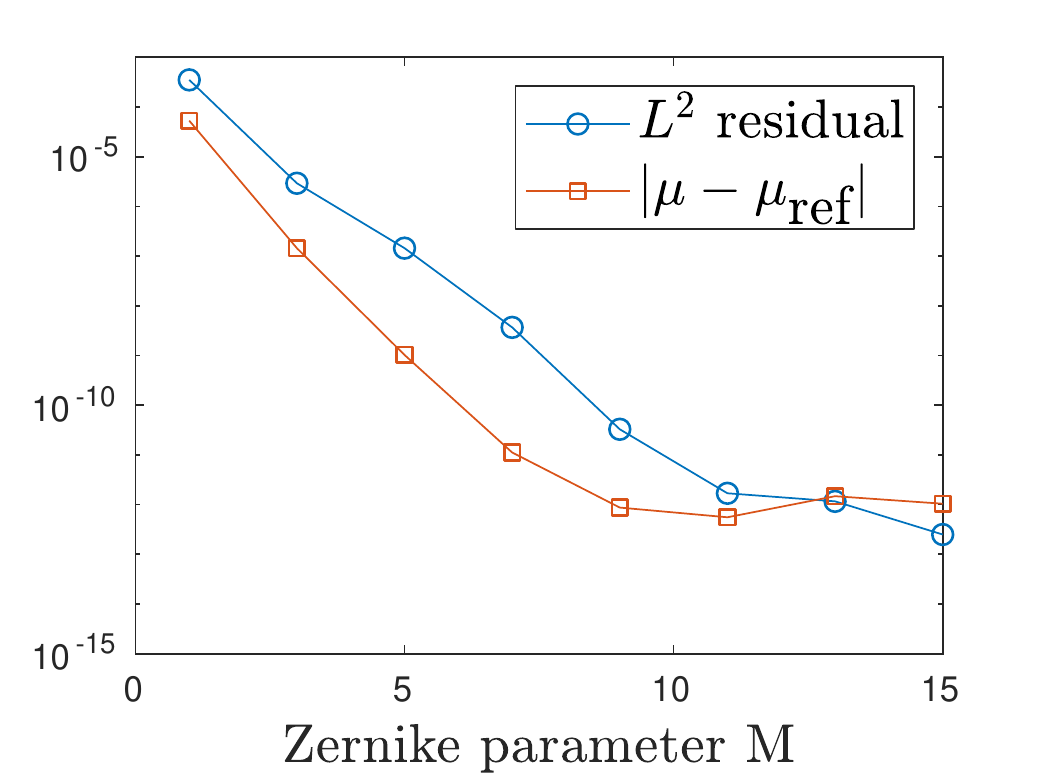}
    \caption{$L^2$ residual error and self-consistent convergence error 
of the average MFPT $\mu$ for random patches with
      $f_{N,\varepsilon} = 0.05$. Left panel: $N = 10$, $\varepsilon
      \approx 0.141$. Middle panel: $N =
    100$, $\varepsilon \approx 0.0447$. Right panel: $N = 1000$,
  $\varepsilon \approx 0.0141$.}
    \label{fig:Mconvergence}
\end{figure}

The results are presented in Fig. \ref{fig:Mconvergence}. In all cases, we
observe the expected spectral convergence with respect to $M$, and 
can reach errors of approximately $10^{-12}$ or less. We
also find that the residual error appears to provide a good upper bound
on the error of $\mu$ until convergence is reached.

\subsection{Large scale simulations} \label{sec:numexamples}

We next study the performance of our solver 
as $N$ is increased and $\varepsilon$ is decreased. 
The error is measured by computing the $L^2$ residual 
\eqref{eq:l2res} on a random
patch. The parameters for the one-patch solver are set as in the 
previous section with $M = 15$, but we fix the ID tolerance at
$10^{-11}$, the GMRES tolerance at $10^{-10}$, and the incoming grid
truncation tolerance at $10^{-8}$. This selection of parameters yields
errors in range $10^{-7}-10^{-10}$ for all of our experiments.
Our calculations are performed on either a laptop with a 4-core Intel
i7-3630QM 2.40GHz processor or a
workstation with four Intel Xeon E7-4880 2.50GHz processors.
each of which has 15 cores.
The algorithm has been implemented in
Fortran, and in both cases, the hierarchical fast algorithm is
parallelized over all available cores using OpenMP.

We consider randomly located patches, uniformly located patches and patches
that are highly clustered.
For each experiment we report $N$, $\varepsilon$, 
the computed value of the average MFPT $\mu$, 
truncated at $8$ significant digits,
the $L^2$ residual error on a random patch,
the total number of GMRES iterations,
the total solve time, and the 
time per GMRES iteration.
We also compute the {\em parallel
scaling factor} - namely, the ratio of the time to compute
the matrix-vector product \eqref{eq:sysmatapply} using a
single core to the time required using all cores on the 60-core
workstation.

\subsubsection{Example 1: Random patches with area fraction
$f_{N,\varepsilon} = 0.05$}

\begin{figure}[ht]
  \centering
    \includegraphics[width=.32\linewidth]{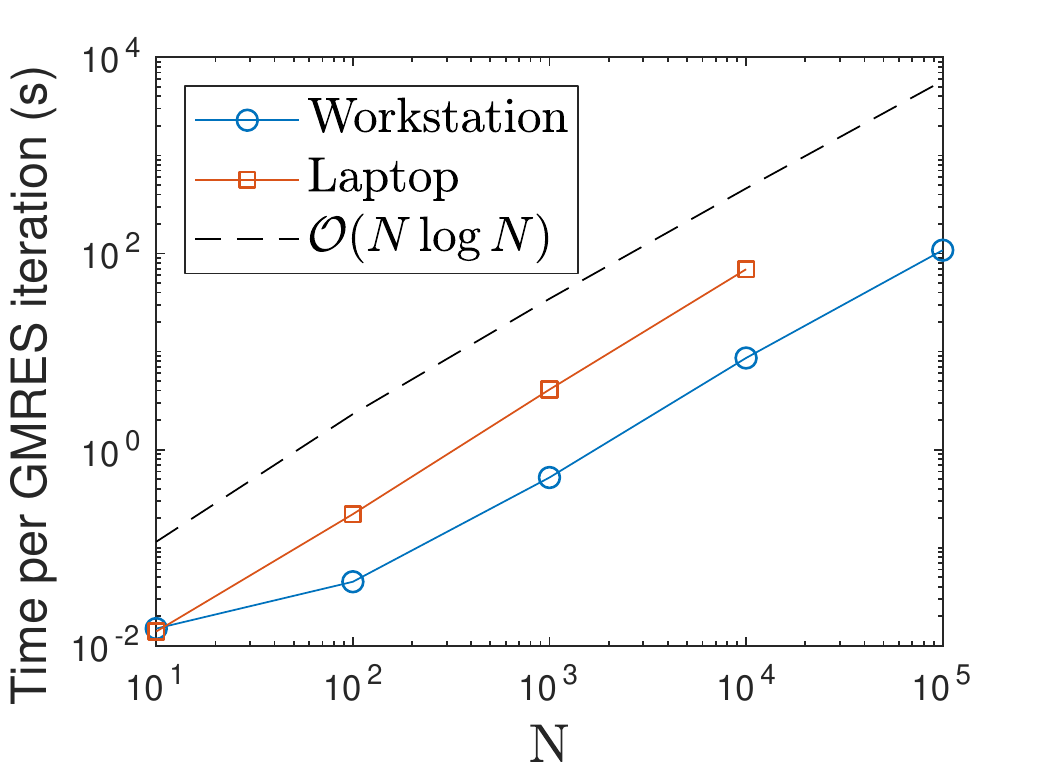} \hspace{.01in}
    \includegraphics[width=.32\linewidth]{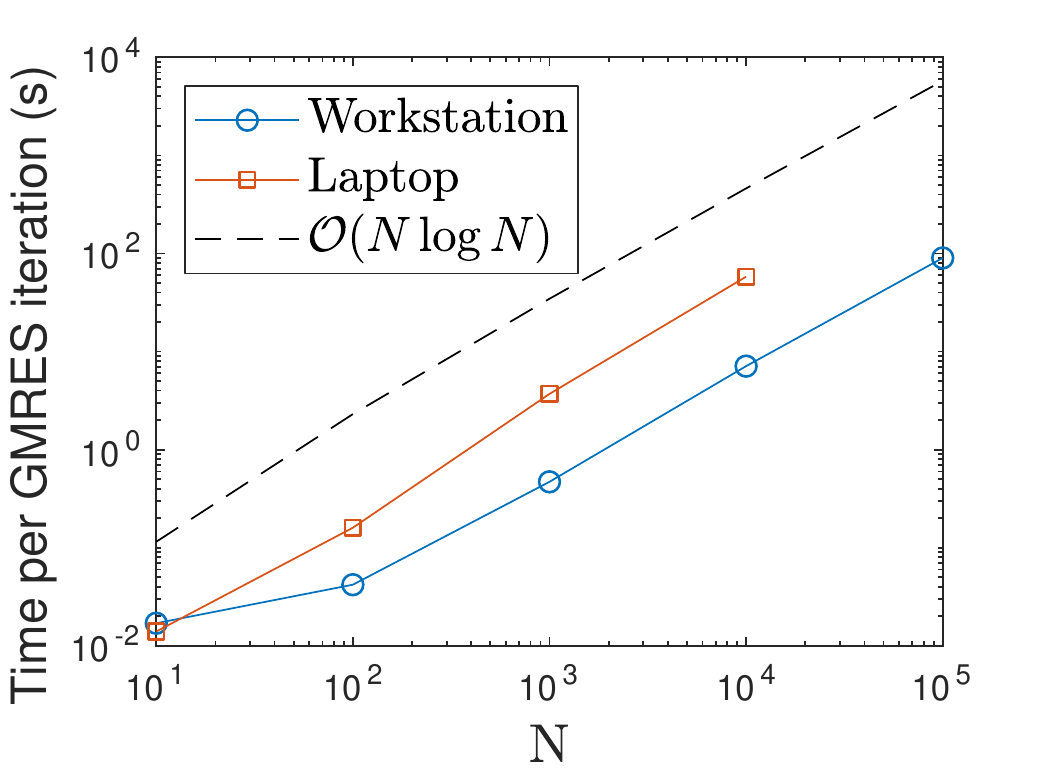} \hspace{.01in}
    \includegraphics[width=.32\linewidth]{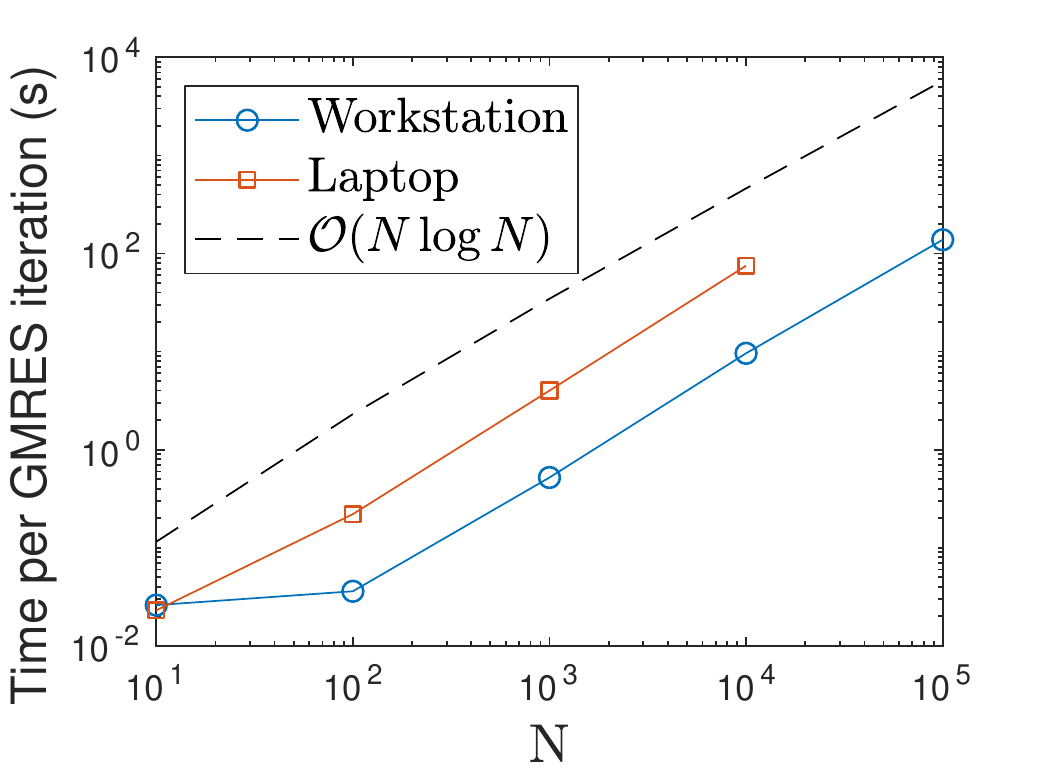}
    Example 1 \hspace{1.4in} Example 2 \hspace{1.4in} Example 3
    \caption{Time per GMRES iteration for the 4-core laptop and 60-core
    workstation. A reference curve with $\OO(N \log N)$ scaling is
    also plotted.} 
    \label{fig:timings}
\end{figure}

Fixing the patch area fraction at $f_{N,\varepsilon} = 0.05$,
we let $\varepsilon$ be given by
\eqref{eq:patareafrac} for $N = 10,100,1000,10\, 000,100\, 000$, 
with patches randomly distributed on the sphere 
with a minimum patch separation of $2
\varepsilon$. 
The corresponding results are given in Table \ref{tab:exrand}. 
In the left panel of Fig. \ref{fig:timings}, we plot the time
per GMRES iteration as a function of $N$ using the 4-core
laptop and the 60-core workstation, as well as a reference curve
with $\OO(N \log N)$ scaling. 
In Fig. \ref{fig:sphereplots}, we also plot the
computed MFPT $\vbar$ just inside the unit sphere - 
on a sphere of radius $1-\varepsilon/5$ - 
for $N= 10,100,1000,10\, 000$.
The case $N = 100\, 000$ case was plotted earlier, 
in Fig. \ref{fig:randpts1e5}.

Note that the number of GMRES iterations increases with $N$,
as one would expect from the increased complexity of the problem,
but slowly.
The computation with $N = 100\, 000$ required just over
an hour to complete using the 60-core workstation. 
The computation with $N = 10\, 000$
required just over $45$ minutes to solve on the 4-core laptop, 
and the computation with $N = 1000$ required approximately
one minute.
(The case $N = 100\, 000$  was not attempted on the laptop because of
memory requirements.) 
Note from the data in Table \ref{tab:exrand} that
we achieve approximately $85\%$ parallel efficiency at 
$N=1000$ and an efficiency near
$90\%$ for the largest calculation.
Note also from Fig. \ref{fig:timings} that the complexity
of the fast algorithm is consistent with the expected 
$O(N \log N)$ scaling. 

\begin{table}[ht]
  \small
  \centering
  \begin{tabular}{|c|r|r|r|r|r|}
    \hline
    $N$ & $10$ & $100$ & $1000$ & $10\, 000$ & $100\, 000$ \\ \hline
    $\varepsilon$ & $\approx 0.14$ & $\approx 0.045$ & $\approx
    0.014$ & $\approx 0.0045$ & $\approx 0.0014$  \\ \hline
    Average MFPT $\mu$ & $0.64277353$ & $0.24999828$ & $0.12308716$ &
    $0.084405945$ & $0.072275200$ \\ \hline
    $L^2$ residual error & $3.6 \times 10^{-9}$ & $1.6 \times 10^{-9}$ & $5.3 \times 10^{-9}$ & $4.8 \times 10^{-8}$ & $2.2 \times 10^{-8}$  \\ \hline
    $\#$ GMRES iterations & $7$ & $12$ & $17$ & $25$ & $35$ \\ \hline
    Total iteration time (s) (60 cores) & $0.11$ & $0.54$ & $8.9$ &
    $215$ & $3793$ \\ \hline
    Time per iteration (s) (60 cores) & $0.02$ & $0.05$ & $0.5$ &
    $8.6$ & $108$  \\ \hline
    Total iteration time (s) (laptop) & $0.10$ & $2.63$ & $68.9$  &
    $1731$ &   \\ \hline
    Time per iteration (s) (laptop)  & $0.01$ & $0.22$ & $4.1$ &
    $69$ &  \\ \hline
    Parallel scaling factor (60 cores) & $2.1$ & $25.7$ & $51.4$ & $52.3$  & $53.5$  \\

    \hline
  \end{tabular}
  \caption{Narrow escape problem with random patches at patch
area fraction $f_{N,\varepsilon}
  = 0.05$.}
  \label{tab:exrand}
\end{table}

\subsubsection{Example 2: Uniform patches with area fraction
$f_{N,\varepsilon} = 0.05$}

Using the same patch area fraction as in the previous example,
we let $N$ take the same values, but place the patch centers at the Fibonacci
spiral points, which are approximately uniform 
on the sphere \cite{Bernoff2018}. 
Results are shown in Table
\ref{tab:exunif} and the middle panel of Fig. \ref{fig:timings}.
The computed MFPT $\vbar$ on the sphere
of radius $1-\varepsilon/5$ was plotted in
Fig.  \ref{fig:fibopts1e4} for the case $N = 10\, 000$. The MFPT is plotted
for the $N = 100$ and $N = 1000$ cases in Fig. \ref{fig:sphereplots}. 

\begin{table}[ht]
  \small
  \centering
  \begin{tabular}{|c|r|r|r|r|r|}
    \hline

    $N$ & $10$ & $100$ & $1000$ & $10\, 000$ & $100\, 000$ \\ \hline
    $\varepsilon$ & $\approx 0.14$ & $\approx 0.045$ & $\approx
    0.014$ & $\approx 0.0045$ & $\approx 0.0014$  \\ \hline
    Average MFPT $\mu$ & $0.62771752$  & $0.23201408$ & $0.11813387$ &
    $0.082870386$ & $0.071784189$ \\ \hline
    $L^2$ residual error & $3.0 \times 10^{-9}$ & $1.5 \times 10^{-9}$ & $3.2 \times 10^{-8}$ & $6.4 \times 10^{-8}$ & $8.4 \times 10^{-8}$   \\ \hline
    $\#$ GMRES iterations & $6$ & $9$ & $11$ & $16$ & $20$ \\ \hline
    Total iteration time (s) (60 cores) & $0.10$ & $0.38$ & $5.1$ &
    $114$ & $1803$ \\ \hline
    Time per iteration (s) (60 cores) & $0.02$ & $0.04$ & $0.47$ &
    $7.1$ & $90$ \\ \hline
    Total iteration time (s) (laptop) & $0.087$ & $1.45$ & $40.7$ &
    $926$ &  \\ \hline
    Time per iteration (s) (laptop) & $0.014$ & $0.16$ & $3.7$ & $58$ &  \\ \hline
    Parallel scaling factor (60 cores) & $5.0$ &
    $29.9$ & $53.7$ & $54.0$ & $54.8$ \\

    \hline
  \end{tabular}
  \caption{Narrow escape problem with uniform patches
   at patch area fraction $f_{N,\varepsilon} = 0.05$.}
  \label{tab:exunif}
\end{table}

\subsubsection{Example 3: Clustered patches}

In our final example, we configure the patches to form a collection 
of $20$ clusters. Each
cluster is contained within a disk on the surface of the sphere centered
at the vertices of a dodecahedron inscribed in
the sphere, and the radii of the disks are chosen so that all $20$ disks
cover one quarter of the area of the sphere. Patch centers are placed
randomly on the sphere, and a proposed center is accepted if it falls
within one of the disks, while enforcing a minimum patch separation distance
of $2 \varepsilon$. We choose $\varepsilon$ empirically to be as large
as possible so that our random placement process yields the desired
number $N$ of patches in a reasonable amount of time. For sufficiently
large $N$, this results in a
much denser packing of patches within each cluster than we had in our
previous examples.

\begin{figure}[p!]
  \centering
    \includegraphics[width=.29\linewidth]{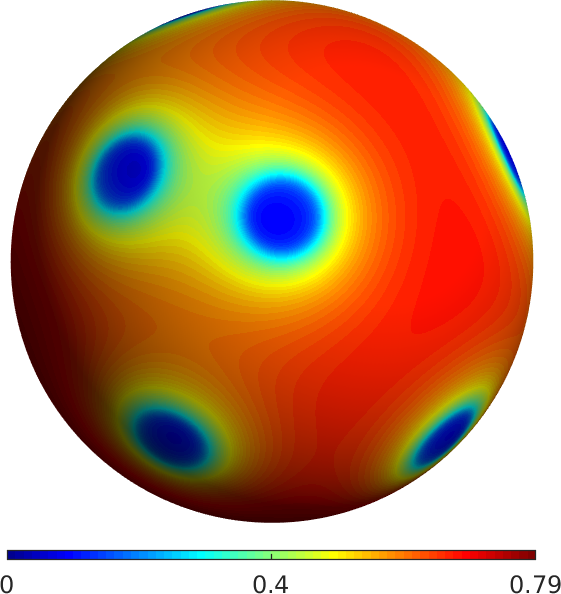} \hspace{.1in}
    \includegraphics[width=.29\linewidth]{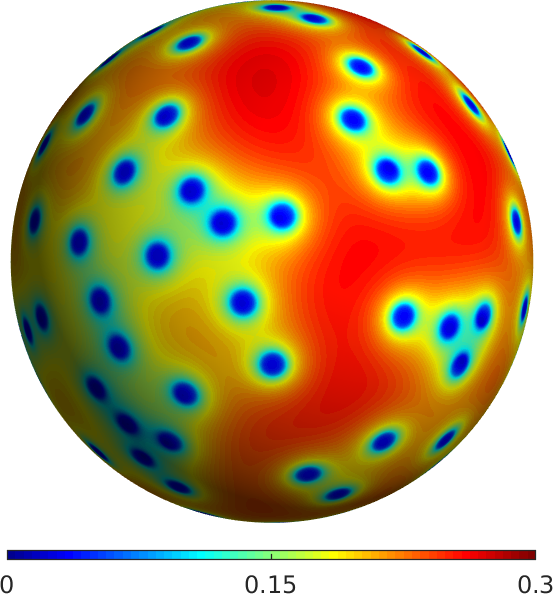} 

    \vspace{.1in}

    \includegraphics[width=.29\linewidth]{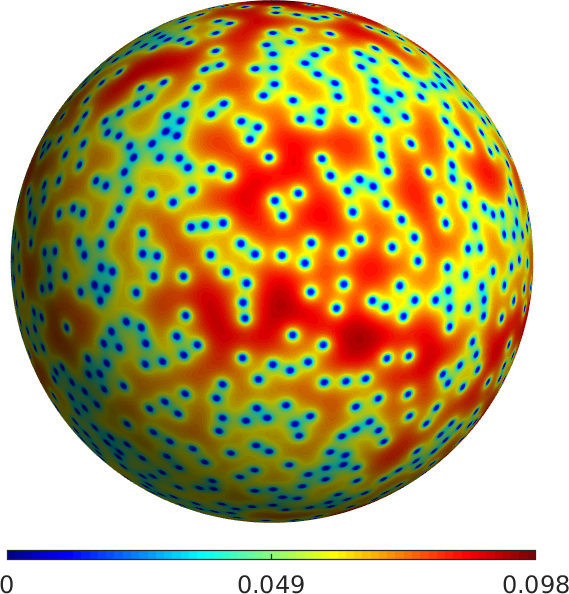} \hspace{.1in}
    \includegraphics[width=.29\linewidth]{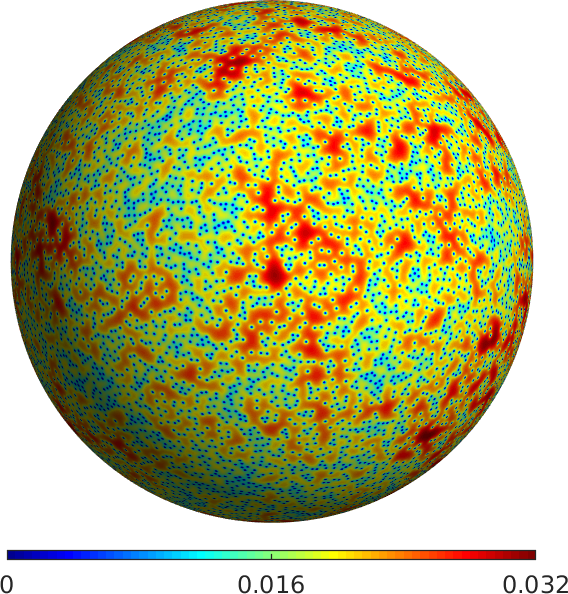}
    
    \vspace{.1in}
    
    \includegraphics[width=.29\linewidth]{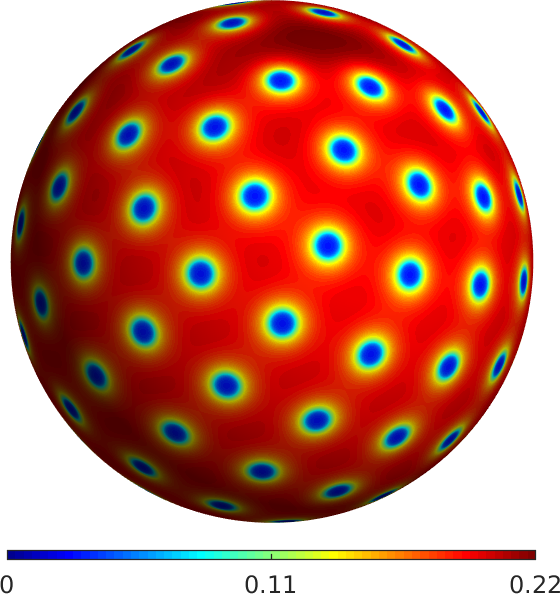} \hspace{.1in}
    \includegraphics[width=.29\linewidth]{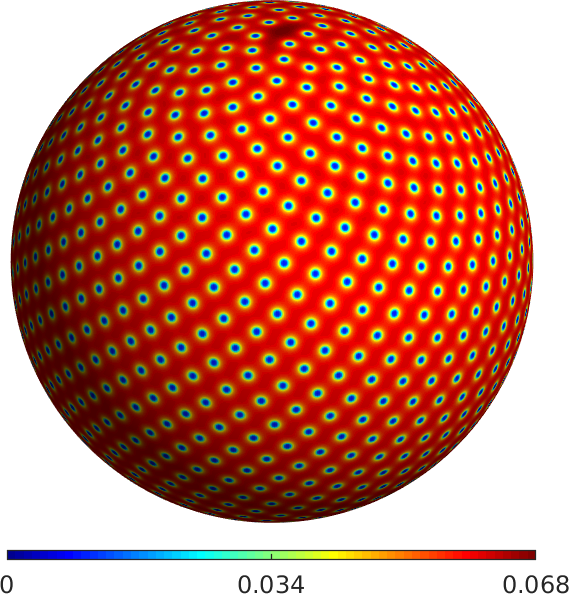}
    
    \vspace{.1in}
    
    \includegraphics[width=.29\linewidth]{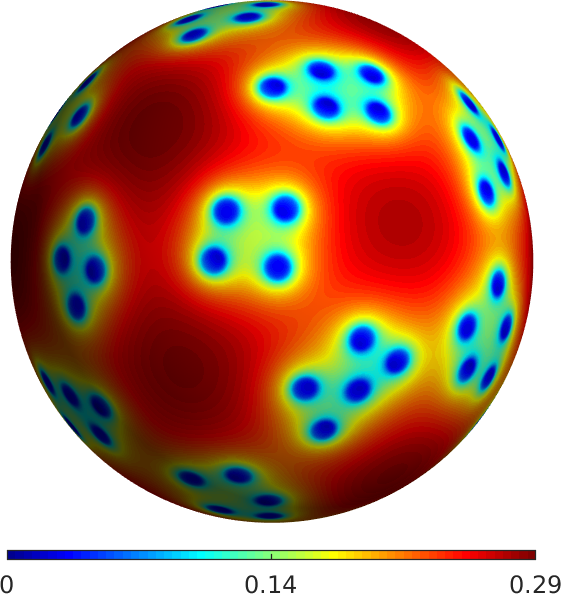} \hspace{.1in}
    \includegraphics[width=.29\linewidth]{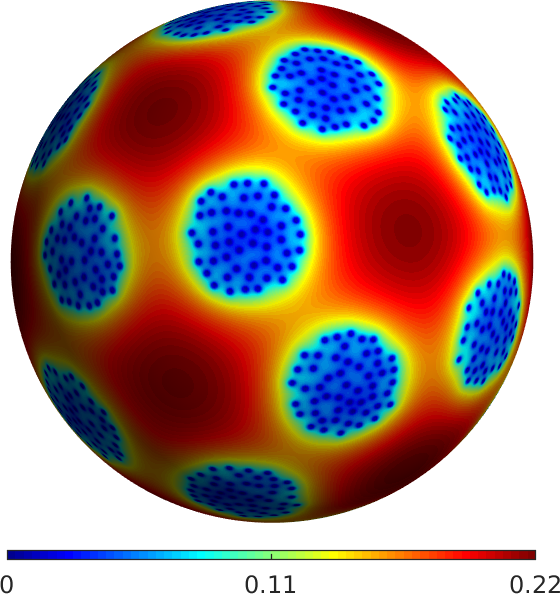}
    
    \caption{Plots of the MFPT $\vbar$ on a sphere of radius $1
      - \varepsilon/5$ for the experiments described in
      Section \ref{sec:numexamples}. The first two rows correspond to
  Example 1 with $N =
  10,100,1000,10\, 000$. The third row corresponds to Example 2 
  with $N = 100, 1000$. The final row corresponds to 
  Example 3 with $N = 100, 1000$.}
    \label{fig:sphereplots}
\end{figure}

The results of our simulations are provided in Table
\ref{tab:exclus} and the right panel of Fig. \ref{fig:timings}.
The MFPT is plotted on a sphere of radius $1-\varepsilon/5$ in Fig. \ref{fig:cluspts1e4} for the $N =
10\, 000$ case and in Fig. \ref{fig:sphereplots} for the $N = 100$ and $N =
1000$ cases. The
denser packing of patches leads to 
a greater number of GMRES iterations than in the previous examples and
longer computation times, but the
difference is mild. The case with $N = 100\, 000$ required 
just over an hour and a half to solve on our 60-core workstation. 
The simulation with $N = 10\, 000$
required 75 minutes on a laptop, and the simulation with $N =
1000$ required about one minute.

\begin{table}[ht]
  \small
  \centering
  \begin{tabular}{|c|r|r|r|r|r|}
    \hline

    $N$ & $10$ & $100$ & $1000$ & $10\, 000$ & $100\, 000$ \\ \hline
    $\varepsilon$ & $0.25$ & $0.047$ & $0.012$ & $0.0035$  & $0.001$ \\ \hline
    Average MFPT $\mu$ & $0.29687267$ & $0.25519357$ & $0.20318506$ &
    $0.17622000$ & $0.16531162$ \\ \hline
    $L^2$ residual error & $4.9 \times 10^{-10}$ & $3.9 \times 10^{-9}$
    & $1.2 \times 10^{-8}$ & $6.5 \times 10^{-8}$ & $1.2 \times 10^{-7}$
    \\ \hline
    $\#$ GMRES iterations & $8$ & $12$ & $19$ & $28$ & $42$ \\ \hline
    Total iteration time (s) (60 cores) & $0.21$ & $0.43$ & $9.9$ &
    $269$ & $5795$   \\ \hline
    Time per iteration (s) (60 cores) & $0.03$ & $0.04$ & $0.52$ &
    $9.6$ & $138$  \\ \hline
    Total iteration time (s) (laptop) & $0.18$ & $2.7$ & $76.4$ &
    $2112$ &  \\ \hline
    Time per iteration (s) (laptop)  & $0.02$ & $0.22$ & $4.0$ & $75$ &  \\ \hline
    Parallel scaling factor (60 cores) & $2.9$ &
    $43.9$ & $49.3$ & $51.4$ & $55.5$ \\

    \hline
  \end{tabular}
  \caption{Narrow escape problem with clustered patches.}
  \label{tab:exclus}
\end{table}

\begin{remark}
We carried out the simulations above for the corresponding 
exterior problem as well (the narrow capture problem).
As expected (since the integral equations are nearly identical),
the timings and errors are similar and are therefore omitted.
\end{remark}

\section{Conclusions}\label{sec:conclusions}

We have developed a fast solver for the narrow capture and narrow escape
problems on the sphere with arbitrarily-distributed well-separated 
disk-shaped patches. We solve the corresponding mixed
boundary value problems by an integral equation scheme
derived using the Neumann Green's functions for the sphere. Our
numerical method combines a high order accurate solver for the one-patch
problem, a multiple scattering formalism, and a hierarchical fast
algorithm. We have demonstrated the scheme on examples with $N$ as large
as $100\, 000$, significantly larger than previously accessible.
The ability to carry out such large-scale simulations will permit
a systematic study of the asymptotic
approaches described, for example, in \cite{Cheviakov2010} and
\cite{Lindsay2017}. 

Possible extensions of our method include the consideration of
narrow escape and narrow capture problems when the 
patches are asymmetric and have multiple shapes.
Assuming some separation between patches, the
multiple scattering formalism still applies, but the
single patch integral
equation will not be solvable by separation of variables and
the compressed representation of outgoing fields will need to be computed
for each distinct patch type. 
Neither of these extra steps, however, affects the asymptotic 
$\OO(N \log N)$ scaling of the fast algorithm.
Exterior problems involving
multiple spheres with different arrangements of patches
could also be simulated by a simple modification of our multiple scattering 
approach. 

A more challenging problem is to extend our method to non-spherical
geometries. For this, one would either have to discretize the entire
domain surface, rather than just the absorbing patches, or construct the
Neumann Green's function for such a domain numerically. 
In the latter case, aspects of
our multiple scattering approach would carry over. We are
currently investigating these issues and will report on our progress
at a later date.

\section*{Acknowledgments}
We would like to thank Michael Ward for suggesting this problem
and for several valuable insights. We would also like to
thank Mike O'Neil for many useful conversations. J.K.
was supported in part by the Research Training Group in Modeling
and Simulation funded by the National Science Foundation via grant
RTG/DMS-1646339.

\bibliographystyle{ieeetr}
{\footnotesize \bibliography{kg_sphere}}

\end{document}